\documentclass[12pt,reqno]{amsart}
\usepackage{amsthm}
\usepackage{amssymb}
\usepackage{graphics}
\usepackage{tikz}
\usepackage{float}
\usetikzlibrary{shapes,backgrounds,calc}
\usepackage{latexsym}
\usepackage{multicol}
\usepackage{verbatim,enumerate}
\usepackage{accents}
\usepackage{cite}
\usepackage{subfig}
\usepackage{multirow}
\usepackage{bigstrut}
\usepackage{amsthm}
\usepackage{amssymb}
\usepackage{graphics}
\usepackage{tikz}
\usetikzlibrary{shapes,backgrounds,calc}
\usepackage{latexsym}
\usepackage{multicol}
\usepackage{verbatim,enumerate}
\usepackage{accents}
\usepackage{cite}
\usepackage{array}
\usepackage[colorlinks=true, linkcolor=blue, citecolor=blue, urlcolor=black]{hyperref}
\usepackage{hyperref}
\usepackage{amsmath, amscd,url}
\usepackage{multirow}
\usepackage{bigstrut}
\usepackage{longtable}

\usepackage{stackengine}

\newcommand{\ind}{\mathop{\mathrm{ind}}}

\newcommand{\F}{\mbox{$\mathbb F$}}	% For Finite fields

\usepackage{hyperref}
\usepackage{amsmath, amscd,url}
\usepackage{longtable}

\advance\textwidth by 1.3in \advance\oddsidemargin by -.6in \advance\evensidemargin by -.6in
\theoremstyle{definition}

\newtheorem{theorem}{Theorem}[section]
\newtheorem{lemma}[theorem]{Lemma}
\newtheorem{corollary}[theorem]{Corollary}
\theoremstyle{definition}
\newtheorem{definition}[theorem]{Definition}
\newtheorem{example}[theorem]{Example}

\newtheorem{prop}{Proposition}[section]
\theoremstyle{remark}

\theoremstyle{definition}

\newcounter{cnt}
 \makeatletter
\def\mydggeometry{\makeatletter\dg@YGRID=1\dg@XGRID=20\unitlength=0.003pt\makeatother}
\makeatother \theoremstyle{remark}

\numberwithin{equation}{section}
\let\bwdg\bigwedge
\def\bigwedge{{\textstyle\bwdg}}

\newcommand{\nc}{\newcommand}
\newcommand{\rnc}{\renewcommand}

\nc{\cal}{\mathcal} \nc{\goth}{\mathfrak} \rnc{\bold}{\mathbf}

\nc\bomega{{\mbox{\boldmath $\omega$}}} \nc\bpsi{{\mbox{\boldmath $\Psi$}}}
\nc\balpha{{\mbox{\boldmath $\alpha$}}}
\nc\bpi{{\mbox{\boldmath $\pi$}}}
\nc\bvpi{{\mbox{\boldmath $\varpi$}}}
\nc\chara{\operatorname{ch}}

\nc\bxi{{\mbox{\boldmath $\xi$}}}
\nc\bmu{{\mbox{\boldmath $\mu$}}} \nc\bcN{{\mbox{\boldmath $\cal{N}$}}} \nc\bcm{{\mbox{\boldmath $\cal{M}$}}} \nc\blambda{{\mbox{\boldmath
			$\lambda$}}}\nc\bnu{{\mbox{\boldmath $\nu$}}}

\makeatletter
\def\section{\def\@secnumfont{\mdseries}\@startsection{section}{1}%
	\z@{.7\linespacing\@plus\linespacing}{.5\linespacing}%
	{\normalfont\scshape\centering}}
\def\subsection{\def\@secnumfont{\bfseries}\@startsection{subsection}{2}%
	{\parindent}{.5\linespacing\@plus.7\linespacing}{-.5em}%
	{\normalfont\bfseries}}
\makeatother

\nc{\Hom}{\operatorname{Hom}}
\nc{\mode}{\operatorname{mod}}
\nc{\End}{\operatorname{End}} \nc{\wh}[1]{\widehat{#1}} \nc{\Ext}{\operatorname{Ext}} \nc{\ch}{\text{ch}} \nc{\ev}{\operatorname{ev}}
\nc{\Ob}{\operatorname{Ob}} \nc{\soc}{\operatorname{soc}} \nc{\rad}{\operatorname{rad}} \nc{\head}{\operatorname{head}}

\nc{\Cal}{\cal} \nc{\Xp}[1]{X^+(#1)} \nc{\Xm}[1]{X^-(#1)}
\nc{\on}{\operatorname} \nc{\Z}{{\bold Z}} \nc{\J}{{\cal J}}  \nc{\Q}{{\bold Q}}

\nc{\N}{{\bold N}}  \nc\boa{\bold a} \nc\bob{\bold b} \nc\boc{\bold c} \nc\bod{\bold d} \nc\boe{\bold e} \nc\bof{\bold f} \nc\bog{\bold g}
\nc\boh{\bold h} \nc\boi{\bold i} \nc\boj{\bold j} \nc\bok{\bold k} \nc\bol{\bold l} \nc\bom{\bold m} \nc\bon{\mathbb n} \nc\boo{\bold o}
\nc\bop{\bold p} \nc\boq{\bold q} \nc\bor{\bold r} \nc\bos{\bold s} \nc\boT{\bold t} \nc\boF{\bold F} \nc\bou{\bold u} \nc\bov{\bold v}
\nc\bow{\bold w} \nc\boz{\bold z}\nc\ba{\bold A} \nc\bb{\bold B} \nc\bc{\mathbb C} \nc\bd{\bold D} \nc\be{\bold E} \nc\bg{\bold
	G} \nc\bh{\bold H} \nc\bi{\bold I} \nc\bj{\bold J} \nc\bk{\bold K} \nc\bl{\bold L} \nc\bm{\bold M} \nc\bn{\mathbb N} \nc\bo{\bold O} \nc\bp{\bold
	P} \nc\bq{\bold Q} \nc\br{\bold R} \nc\bs{\bold S} \nc\bt{\bold T} \nc\bu{\bold U} \nc\bv{\bold V} \nc\bw{\bold W} \nc\bz{\mathbb Z} \nc\bx{\bold
	x} \nc\KR{\bold{KR}} \nc\rk{\bold{rk}} \nc\het{\text{ht }}

\nc\toa{\tilde a} \nc\tob{\tilde b} \nc\toc{\tilde c} \nc\tod{\tilde d} \nc\toe{\tilde e} \nc\tof{\tilde f} \nc\tog{\tilde g} \nc\toh{\tilde h}
\nc\toi{\tilde i} \nc\toj{\tilde j} \nc\tok{\tilde k} \nc\tol{\tilde l} \nc\tom{\tilde m} \nc\ton{\tilde n} \nc\too{\tilde o} \nc\toq{\tilde q}
\nc\tor{\tilde r} \nc\tos{\tilde s} \nc\toT{\tilde t} \nc\tou{\tilde u} \nc\tov{\tilde v} \nc\tow{\tilde w} \nc\toz{\tilde z} \nc\woi{w_{\omega_i}}

\begin{document}
	\setcounter{section}{0}
	\setcounter{tocdepth}{1}

	%%%%%%%%%%%%%%%%%%%%%%%%%%%%%%%%%%%%%%%%%%%%

\title{DISCRIMINANT and integral basis OF $\Q(\sqrt[12]{a})$}
\author[Surender Kumar]{Surender Kumar}
\author[Anuj Jakhar]{Anuj Jakhar}
%\address[Anuj Jakhar]{Indian Institute of Technology (IIT) Bhilai}
%\email[Anuj Jakhar]{anujjakhar@iitbhilai.ac.in \\ anujiisermohali@gmail.com}

\thanks{The first author is grateful to the University Grants Commission, New Delhi for providing financial support in the form of Junior Research Fellowship through Ref No.1129/(CSIR-NET JUNE 2019). The second author is thankful to SERB grant SRG/2021/000393.}


\begin{abstract}
	\noindent  Suppose $m$ be a $12$-th power free integer. Let $K=\Q(\theta)$ be an  algebraic number field defined by a complex root $\theta$ of an irreducible polynomial $x^{12}-m$ and $O_K$ be its ring of integers. In this paper,  we determine the highest power of $p$ dividing the index of the subgroup $\Z[\theta]$ in $O_K$ and   $p$-integral basis of $K$ for each prime $p$. These $p$-integral bases lead to the construction of an
	integral basis of $K$ which is illustrated with examples. In particular, when $m$ is a square free integer, we provide necessary and sufficient conditions for the set $\{1,\theta,\theta^2,\cdots,\theta^{10},\theta^{11}\}$ to be an integral basis of $K.$ 
\end{abstract}
\maketitle

\section{Introduction}\label{intro}
Let $f(x)\in\Z[x] $ be a monic irreducible polynomial of degree $n$ over the field  $\Q$ of rationals. Let $K=\Q(\theta)$ be an algebraic number field with $\theta$ a root of $f(x)$ and $O_K$ be its ring of algebraic integers. The computation of discriminant and construction of an integral basis of a number field are the major problems in algebraic number theory. Recently, many mathematicians (\cite{1}, \cite{3}, \cite{4},\cite{a}, \cite{11}, \cite{12})  showed their interest in this research area.
 %We denote  the index of the subgroup $\Z[\alpha]$ in $O_K$  by $\ind\alpha,$ for any algebraic integer $\alpha\in K$. The $\Z$-basis of $O_K$ is called integral basis of $K$. 
 %If  for some $\gamma\in O_K$, $\ind\gamma=1$ then  $\{1,\gamma,\gamma^2,\cdots,\gamma^{n-1}\}$ is an integral basis of $K$. In such case, we say  the field $K$ is monogenic.   By a pure algebraic number field of degree $r$, we mean an algebraic number field of the form $\Q(\sqrt[r]{a})$ defined by a root  of an irreducible polynomial $x^r-a\in\Z[x].$
  In 1897, Landsberg \cite{2} gave a formula for the discriminant of pure prime degree number fields. In 1900, Dedekind \cite{1} gave an explicit integral basis for pure cubic fields. In $1984$, Funakura \cite{3}  provided an integral basis and  a formula for the discriminant of all pure quartic fields. In $2015$,
Hameed and Nakahara \cite{4} gave a formula for the discriminant of pure octic number fields $\Q(\sqrt[8]{a})$, where $a$ is a squarefree integer.  In 2022, Jakhar \cite{9} provided an explicit $p$-integral basis for the number field defined by a root of an irreducible polynomial $x^{p_1p_2}-a\in\Z[x],$ where $p_1$ and $p_2$ are distinct  primes. \\
\indent Let $p$ be a prime number and  $\Z_{(p)}$ denote the localization of the ring $\Z$ of integers at the
prime ideal $p\Z.$ Then the integral closure $R_{(p)}$ of $\Z_{(p)}$ in an algebraic number field $K$ of degree $n$ is a free module over 
$\Z_{(p)}$ having rank $n$. A $\Z_{(p)}$-basis of the module
$R_{(p)}$ is called a $p$-integral basis of $K$.\\
 %These $p$-integral bases lead to a construction of an explicit integral basis of $K$ taking into considerations all the primes $p$ dividing the index of the subgroup $\Z[\theta]$ in $O_K$.\\
  %\indent In $2017$, Jakhar, Khanduja and Sangwan, \cite{6} gave a formula for the discriminant of all those $n$-th degree pure number fields $\Q(\sqrt[n]{a})$ which are such that for each prime $p$ dividing $n$, either $p\nmid a$ or $p$ does not divide $v_p(a)$, where $v_p(a)$ stands for the highest power of p dividing a; clearly this condition is satisfied when either $a$; $n$ are coprime or $a$ is squarefree. But this result does not provide any conclusion, when a rational prime $p$, divides $n$ and $\gcd(p,v_p(a))$. Recall that  the discriminant of $g(x)$ and the discriminant $d_L$ of $L$ are related by the formula
  %\begin{equation}\label{eq1}
  %	discr(g)=(\ind \xi)^2 d_L. 
  %\end{equation}
   %for the algebraic number field $L=\Q(\xi)$ with $\xi$ an algebraic integer having minimal polynomial $g(x)$ over  $\Q.$
  \indent In the present paper, let  $K=\Q(\theta)$ be an algebraic number field with $\theta$ a root of an irreducible polynomial $f(x)=x^{12}-m,$ where $m$ is a $12$-th power free integer. Our goal is to determine the highest power of any  prime $p$ dividing the index of the subgroup $\Z[\theta]$ in $O_K$ and  to construct a $p$-integral basis of $K$. We shall denote the index of the subgroup $\Z[\theta]$ in $O_K$ by $\ind\theta.$ It is easy to check that  the discriminant $D_f$ of the polynomial $f(x)$ is given by 
\begin{equation}\label{eq2}
	D_f=-2^{24} 3^{12}m^{11}.
	\end{equation}
\section{Main results and Examples} 
\indent In what follows, for a prime $p$ and a nonzero
$t$ belonging to the ring $\Z_p$ of $p$-adic integers, $v_p(t)$ will denote the highest power of
$p$ dividing $t$. Let $c_p= \frac{c}{p^{v_p(c)}},$ for any non-zero integer $c$ and rational prime $p$. For a  prime $p$ and  integers $a, b$, we denote $a\equiv b \mod p$ by $a\equiv b~(p).$
  \begin{theorem}\label{Th1.2}
  	Let $K=\Q(\theta)$ be an algebraic number field with $\theta$ a root of an irreducible polynomial $f(x)=x^{12}-m$, where $m$ is a $12$-th power free integer. Let $ p\in\{2,3\}$ be a prime number dividing  both $m$ and $v_p(m)$.  Then $v_p(\ind\theta)$ and a $p$-integral basis are given in Table $1$ and Table $2$ for the primes $p=2$ and $p=3$, respectively. 
  	\vspace{-0.05in}
\begin{center}
	\begin{longtable}[h!]{ m{.7cm} m{.85cm} m{2cm} m{1.33cm} m{9.2cm}}
		\caption{$2$-integral basis with value $v_2(\ind\theta).$}\\ 
		\hline
		\hspace{-0.08in}Case &\hspace{-0.09in}$v_2(m)$& Conditions &$v_2(\ind\theta)$&\hspace{2.5cm} $2$-integral basis\\
		\hline
		
		A1&$2$&$m_2\equiv 1~(4)$& $12$& \tiny \{ $1,$ $\theta,$ $\theta^2,$ $\theta^3,$ $\theta^4,$ $\theta^5,$ $\frac{\theta^6-2}{2^2},$ $\frac{\theta^7-2\theta}{2^2},$ $\frac{\theta^8-2\theta^2}{2^2},$ $\frac{\theta^9-2\theta^3}{2^2},$ $\frac{\theta^{10}-2\theta^4}{2^2},$ $ \frac{\theta^{11}-2\theta^5}{2^2}$\}\\
			
		A2 &$2$&$m_2\equiv3~(8)$ & $13$&\tiny \{ $1$, $\theta,$ $\theta^2,$ $\theta^3,$ $\theta^4,$ $\theta^5,$ $\frac{\theta^6-2\theta^3+6}{2^2},$ $\frac{ \theta^7-2\theta^4+6\theta}{2^2},$ $\frac{\theta^8-2\theta^5+6\theta^2}{2^2},$ $\frac{\theta^9-2\theta^6+6\theta^3}{2^2},$
		$\frac{\theta^{10}-2\theta^7+6\theta^4}{2^2},$  $\frac{\theta^{11}-2\theta^8+6\theta^5}{2^3}$\}\\
		
		A3 &$2$&$m_2\equiv7~(8)$ & $15$&\tiny \{ $1$, $\theta,$ $\theta^2,$ $\theta^3,$ $\theta^4,$ $\theta^5,$ $\frac{\theta^6-2\theta^3+2}{2^2},$ $\frac{ \theta^7-2\theta^4+2\theta}{2^2},$ $\frac{\theta^8-2\theta^5+2\theta^2}{2^2},$ $\frac{\theta^9-2\theta^6+2\theta^3}{2^3},$
		$\frac{\theta^{10}-2\theta^7+2\theta^4}{2^3},$ $ \frac{\theta^{11}-2\theta^8+2\theta^5}{2^3}$\}\\
		
		A4&$4$&$m_2\equiv 3~(4)$& $21$&\tiny\{$1$, $\theta,$ $\theta^2,$ 
		$\frac{\theta^3}{2},$ $\frac{\theta^4}{2},$
		$\frac{\theta^5}{2},$ $\frac{\theta^6+4\theta^3+12}{2^2},$ $\frac{\theta^7+4\theta^4+12\theta}{2^2},$ $\frac{\theta^8+4\theta^5+12\theta^2}{2^3},$ $ \frac{\theta^9+2\theta^6+4\theta^3+8}{2^3},$ $\frac{\theta^{10}+2\theta^7+4\theta^4+8\theta}{2^4},$
		 $\frac{\theta^{11}+2\theta^8+4\theta^5+8\theta^2}{2^4}$\}\\
		
		%\cline{1-1}
		%\cline{3-5}
		A5 &$4$& $m_2\equiv5~(8)$ &$26$& \tiny \{$1$, $\theta,$ $\theta^2,$ 
		$\frac{\theta^3+6}{2},$ $\frac{\theta^4+6\theta}{2},$ 
		$\frac{\theta^5+6\theta^2}{2^2},$ $\frac{\theta^6+4\theta^3+12}{2^3},$ $\frac{\theta^7+4\theta^4+12\theta}{2^3},$  $\frac{\theta^8+4\theta^5+12\theta^2}{2^3},$ $\frac{\theta^9+2\theta^6+4\theta^3+8}{2^4},$ $\frac{\theta^{10}+2\theta^7+4\theta^4+8\theta}{2^4},$ $\frac{\theta^{11}+2\theta^8+4\theta^5+8\theta^2}{2^5}$\}\\
		%\cline{1-1}
		%\cline{3-5}
		A6 &$4$& $m_2\equiv1~(8)$ & $28$&\tiny \{$1$, $\theta, \theta^2,$
		$\frac{\theta^3+6}{2}$, $\frac{\theta^4+6\theta}{2},$
		$\frac{\theta^5+6\theta^2}{2^2},$  $\frac{\theta^6+4\theta^3+12}{2^3},$ $\frac{\theta^7+4\theta^4+12\theta}{2^3},$ 
		$\frac{\theta^8+4\theta^5+12\theta^2}{2^3},$  $\frac{\theta^9+2\theta^6+4\theta^3+8}{2^5},$	
		  $\frac{\theta^{10}+2\theta^7+4\theta^4+8\theta}{2^5},$ $\frac{\theta^{11}+2\theta^8+4\theta^5+8\theta^2}{2^5}$\}\\
		%\cline{1-1}
		%\cline{3-4}
		
		A7&$6$&$m_2\equiv 1~(4)$& $36$& \tiny \{$
		1, \theta, \frac{\theta^2}{2},
		\frac{\theta^3}{2}, \frac{\theta^4}{2^2}, \frac{\theta^5}{2^2}, ~\frac{\theta^6+8}{2^4}, ~\frac{\theta^7+8\theta}{2^4}, \frac{ q'_1(\theta)}{2^5},\frac{\theta q'_1(\theta)}{2^5},    
		\frac{q_1(\theta)}{2^6}, \frac{\theta q_1(\theta)}{2^6}$\},  $q_1(\theta)=\theta^{10}-2\theta^9+2\theta^8-4\theta^6+8\theta^5-8\theta^4+16\theta^2+32\theta+32$ and  $q'_1(\theta)=\theta^{8}-2\theta^7+2\theta^6-4\theta^5+8\theta^4-8\theta^3+8\theta^2-16\theta+16$\\
		%\cline{1-1}
		%\cline{3-4}
		A8&$6$&$m_2\equiv 3~(8)$& $36$&\tiny \{$ 1,\theta,\frac{\theta^2}{2},
		\frac{\theta^3}{2}, \frac{\theta^4}{2^2}, \frac{\theta^5}{2^2}, \frac{s(\theta)}{2^4}, \frac{r(\theta)}{2^4}, \frac{ q'_1(\theta)}{2^5},\frac{\theta q'_1(\theta)}{2^5},    
		\frac{q_1(\theta)}{2^6}, \frac{\theta q_1(\theta)}{2^6}$\},
		 $q_1(\theta)=\theta^{10}+2\theta^8+4\theta^6+8\theta^4+16\theta^2+32$, $q_1'(\theta)=\theta^8-2\theta^6+8\theta^2-16$, $s(\theta)=\theta^6-4\theta^5+8\theta^4+4\theta^3-8$ and $r(\theta)=\theta^7+\theta^6-4\theta^5+4\theta^4+4\theta^3+8\theta-8$ \\
		%\cline{1-1}
		%\cline{3-5}
		%A8 &$6$& $m_2\equiv3~(8)$ &$36$& \tiny$
		%\frac{q'_2(\theta)}{2^2},$ $\frac{\theta q'_2(\theta)}{2^2},$
	%	$\frac{\theta^2 q'_2(\theta)}{2^3},$  $\frac{\theta^3 q'_2(\theta)}{2^3},$  $\frac{ q_2(\theta)}{2^4},$  $\frac{\theta q_2(\theta)}{2^4},$ $\frac{q'_1(\theta)}{2^5},$ $\frac{\theta q'_1(\theta)}{2^5},$  $\frac{q_1(\theta)}{2^6},$  
	%	$ \frac{\theta q_1(\theta)}{2^6},$  $\frac{\theta^2 q'_1(\theta)}{2^6},$ $\frac{\theta^3 q'_1(\theta)}{2^6}$, $q_1(\theta)=\theta^{10}+2\theta^8+4\theta^6+8\theta^4+16\theta^2+32$, $q_2(\theta)=\theta^8+4\theta^6+12\theta^4+32\theta^2+80,$  $q_1'(\theta)=\theta^8-2\theta^6+8\theta^2-16$ and $q'_2(\theta)=\theta^4-4\theta^2+4$\\
		
		A9 &$6$& $m_2\equiv7~(8)$ &$39$&\tiny \{$ 1,\theta,\frac{\theta^2}{2}, \frac{\theta^3}{2}, \frac{\theta^4}{2^2}, \frac{\theta^5}{2^2}, \frac{q_4(\theta)}{2^4}, \frac{q_3(\theta)}{2^5}, \frac{q_1'(\theta)}{2^5}, \frac{\theta q_1'(\theta)}{2^6},\frac{q_1(\theta)}{2^6}, \frac{\theta q_1(\theta)}{2^7}$\}, $q_1(\theta)=\theta^{10}-2\theta^9+2\theta^8-4\theta^6+8\theta^5-8\theta^4+16\theta^2+32\theta+32$,  $q'_1(\theta)=\theta^{8}-2\theta^7+2\theta^6-4\theta^5+8\theta^4-8\theta^3+8\theta^2-16\theta+16$, $ q_3(\theta)=\theta^7-2\theta^6+4\theta^4-8\theta^3-24\theta+16,$ and $q_4(\theta)=\theta^6-4\theta^5+8\theta^4-12\theta^3-16\theta+8$\\
		%\cline{1-1}
		%\cline{3-4}
		
		A10&$8$&$m_2\equiv1~(8) $& $50$&\tiny \{$1$, $\theta,$  $\frac{\theta^2}{2},$
		$\frac{\theta^3-12}{2^2},$ $\frac{\theta^4-12\theta}{2^3},$
		$\frac{\theta^5-12\theta^2}{2^3},$  $\frac{\theta^6-8\theta^3+48}{2^5},$ $\frac{\theta^7-8\theta^4+48\theta}{2^5},$ $\frac{\theta^8-8\theta^5+48\theta^2}{2^6},$ $\frac{\theta^9-4\theta^6+16\theta^3+8}{2^8},$ $\frac{\theta^{10}-4\theta^7+16\theta^4+8\theta}{2^8},$ $\frac{\theta^{11}-4\theta^8+16\theta^5+8\theta^2}{2^9}$\}\\
		%\cline{1-1}
		%\cline{3-5}
		A11 &$8$& $m_2\equiv5~(8)$ & $48$&\tiny \{$1, \theta, \frac{\theta^2}{2},  
		\frac{\theta^3-12}{2^2}, \frac{\theta^4-12\theta}{2^3}, 
		\frac{\theta^5-12\theta^2}{2^3},   \frac{\theta^6-8\theta^3+48}{2^5}, \frac{\theta^7-8\theta^4+48\theta}{2^5},
		\frac{\theta^8-8\theta^5+48\theta^2}{2^6}, $ $\frac{\theta^9-4\theta^6+16\theta^3-64}{2^7}, \frac{\theta^{10}-4\theta^7+16\theta^4-64\theta}{2^8},$ $\frac{\theta^{11}-4\theta^8+16\theta^5-64\theta^2}{2^8}$\}\\
		%\cline{1-1}
		%\cline{3-5}
		A12 &$8$& $m_2\equiv 3~(4)$& $43$&\tiny \{$1$, $\theta,$  $\frac{\theta^2}{2},$ 
		$\frac{q_3(\theta)}{2^2},$ $\frac{\theta q_3(\theta)}{2^2},$
		$\frac{\theta^2 q_3(\theta)}{2^3},$  $\frac{q_2(\theta)}{2^4},$ $\frac{\theta q_2(\theta)}{2^5},$
		$\frac{\theta^2 q_2(\theta)}{2^5},$  $\frac{q_1(\theta)}{2^6},$ $ \frac{\theta q_1(\theta)}{2^7},$ $\frac{\theta^2q_1(\theta)}{2^8}$\}, $q_1(\theta)=\theta^9-4\theta^6+16\theta^3-64$, $q_2(\theta)=\theta^6-8\theta^3+48,$ and $q_3(\theta)=\theta^3-12$\\
		%\cline{1-1}
		%\cline{3-5}
		
		A13&$10$&$m_2\equiv 1~(4)$& $56$&\tiny \{$1$, $\theta,$ $\frac{\theta^2}{2},$ $\frac{\theta^3}{2^2},$ $\frac{\theta^4}{2^3},$ $\frac{\theta^5}{2^4},$ $\frac{\theta^6-32}{2^6},$ $\frac{\theta^7-32\theta}{2^6},$ $\frac{\theta^8-32\theta^2}{2^7},$ $\frac{\theta^9-32\theta^3}{2^8},$ $\frac{\theta^{10}-32\theta^4}{2^9},$ $ \frac{\theta^{11}-32\theta^5}{2^{10}}$\}\\
		%\cline{1-1}
		%\cline{3-5}
		A14 &$10$& $m_2\equiv3~(8)$ & $57$&\tiny \{$1$, $\theta,$ $\frac{\theta^2}{2},$ $\frac{\theta^3}{2^2},$ $\frac{\theta^4}{2^3},$ $\frac{\theta^5}{2^4},$ $\frac{\theta^6-8\theta^3+96}{2^6},$ $\frac{\theta^7-8\theta^4+96\theta}{2^{7}},$ $\frac{\theta^8-8\theta^5+96\theta^2}{2^{7}},$ $\frac{\theta^9-8\theta^6+96\theta^3}{2^{8}},$ $\frac{\theta^{10}-8\theta^7+96\theta^4}{2^{9}},$ $ \frac{\theta^{11}-8\theta^8+96\theta^5}{2^{{10}}}$\}\\
		%\cline{1-1}
		%\cline{3-5}
		A15 &$10$& $m_2\equiv7~(8)$ & $59$&\tiny \{$1$, $\theta,$ $\frac{\theta^2}{2},$ $\frac{\theta^3}{2^2},$ $\frac{\theta^4}{2^3},$ $\frac{\theta^5}{2^4},$ $\frac{\theta^6-8\theta^3+32}{2^6},$ $\frac{\theta^7-8\theta^4+32\theta}{2^{7}},$ $\frac{\theta^8-8\theta^5+32\theta^2}{2^{8}},$ $\frac{\theta^9-8\theta^6+32\theta^3}{2^{9}},$ $\frac{\theta^{10}-8\theta^7+32\theta^4}{2^{9}},$ $ \frac{\theta^{11}-8\theta^8+32\theta^5}{2^{{10}}}$\}\\
		\hline
			\end{longtable}
\end{center}
\vspace{-0.7cm}
\begin{center}
	\begin{longtable}[h!]{ m{.7cm} m{.85cm}m{2cm}m{1.33cm} m{9.2cm}} 
		\caption{$3$-integral basis with value $v_3(\ind\theta).$}\\
		\hline
		\hspace{-0.08in}Case & \hspace{-0.09in}$v_3(m)$&\hspace{-0.05in}Conditions& $v_3(\ind\theta)$&\hspace{2cm} $3$-integral basis\\
		\hline
		B1&$3$ &$m_3\equiv\delta~(9)$& $18$&\tiny \{$1$, $\theta,$ $\theta^2,$ $\theta^3,$ $\frac{\theta^4+6\delta}{3},$ $\frac{\theta^5+6\delta\theta}{3},$ $\frac{\theta^6+6\delta\theta^2}{3^2},$ $\frac{\theta^7+6\delta\theta^3}{3^2},$ $\frac{\theta^8+3\delta\theta^4+9}{3^3},$ $ \frac{\theta^9+3\delta\theta^5+9\theta}{3^{3}},$ $\frac{\theta^{10}+3\delta\theta^6+9\theta^2}{3^3},$ $\frac{\theta^{11}+3\delta\theta^7+9\theta^3}{3^{3}}$\}, $\delta\in\{-1,1\}$\\	
		%\cline{1-1}
		%\cline{3-5}
		B2 &3& $m_3\not\equiv \delta ~(9)$ & $15$&\tiny \{$1$, $\theta,$ $\theta^2,$ $\theta^3,$ $\frac{\theta^4+6\delta}{3},$ $\frac{\theta^5+6\delta\theta}{3},$ $\frac{\theta^6+6\delta\theta^2}{3},$ $\frac{\theta^7+6\delta\theta^3}{3^2},$ $\frac{\theta^8+3\delta\theta^4+9}{3^2},$ $ \frac{\theta^9+3\delta\theta^5+9\theta}{3^2},$ $\frac{\theta^{10}+3\delta\theta^6+9\theta^2}{3^3},$ $\frac{\theta^{11}+3\delta\theta^7+9\theta^3}{3^{3}}$\}, $\delta\in\{-1,1\}.$\\	
		
		B3&$6$ &$m_3\equiv 1~(9)$& $36$&\tiny \{$1,\theta, \frac{\theta^2}{3}, \frac{\theta^3}{3},\frac{\theta^4}{3^2}, \frac{\theta^5-9\theta}{3^3}, \frac{ \theta^6-9\theta^2}{3^3}, \frac{ z_1(\theta)}{3^4},\frac{z_2(\theta)}{3^5}, \frac{\theta z_2(\theta)}{3^5}, \frac{q_1(\theta)}{3^6}, \frac{\theta q_1(\theta)}{3^6}$\}, $q_1(\theta)=\theta^{10}+3\theta^8+9\theta^6+27\theta^4+81\theta^2+243$, $z_1(\theta)=\theta^7-15\theta^5-9\theta^3-27\theta$ and $z_2(\theta)=\theta^8-18\theta^4-162$\\
		%\cline{1-1}
		%\cline{3-5}
		B4 &$6$& $m_3\equiv4~(9)$ & $32$&\tiny \{$ 1,\theta, \frac{\theta^2}{3}, \frac{\theta^2}{3},\frac{\theta^4}{3^2}, \frac{\theta^5}{3^2}, \frac{\theta^6}{3^3}, \frac{\theta^7}{3^3},\frac{\theta^8}{3^4}, \frac{\theta^9-18\theta^5-162\theta}{3^5}, \frac{q_1(\theta)}{3^5}, \frac{\theta q_1(\theta)}{3^6}$\},\\ 
		&&$m_3\equiv 7~(9)$&&\tiny $q_1(\theta)=\theta^{10}+3\theta^8+9\theta^6+27\theta^4+81\theta^2+243$\\
		%\cline{1-1}
		%\cline{3-5}
		B5 &$6$& $m_3\equiv8~(9)$ & $36$& \tiny \{$1$, $\theta,$ $\frac{\theta^2}{3},$ $\frac{\theta^3}{3},$  $\frac{\theta^4-18}{3^2},$ $\frac{\theta^5-18\theta}{3^3},$ $\frac{\theta^6-18\theta^2}{3^3},$ $\frac{\theta^7-18\theta^3}{3^4},$ $\frac{\theta^8-9\theta^4+81}{3^5},$ $\frac{\theta^9-9\theta^5+81\theta}{3^5},$ $\frac{\theta^{10}-9\theta^5+81\theta^2}{3^6},$ $\frac{\theta^{11}-9\theta^6+81\theta^3}{3^6}$\}\\
		%\cline{1-1}
		%\cline{3-5}
		B6 &$6$& $m_3\equiv2~(9)$ & $32$& \tiny \{$1$, $\theta,$ $\frac{\theta^2}{3},$ $\frac{\theta^3}{3},$  $\frac{\theta^4-18}{3^2},$ $\frac{\theta^5-18\theta}{3^2},$ $\frac{\theta^6-18\theta^2}{3^3},$ $\frac{\theta^7-18\theta^3}{3^3},$ $\frac{\theta^8-9\theta^4+81}{3^4},$ $\frac{\theta^9-9\theta^5+81\theta}{3^5},$ $\frac{\theta^{10}-9\theta^5+81\theta^2}{3^5},$ $\frac{\theta^{11}-9\theta^6+81\theta^3}{3^6}$\}\\ &&$m_3\equiv5~(9)$&&\\
		%\cline{1-1}
		%\cline{3-4}
		%\hline
		B7 &$9$ &$m_3\equiv \delta~(9)$& $51$&\tiny \{$1$, $\theta,$ $\frac{\theta^2}{3}, $ $\frac{\theta^3}{3^2},$ $\frac{\theta^4-54\delta}{3^3},$ $\frac{\theta^5-54\delta\theta}{3^4},$ $\frac{\theta^6-54\delta\theta^2}{3^5},$ $\frac{\theta^7-54\delta\theta^3}{3^5},$ $\frac{\theta^8+27\delta\theta^4+729}{3^7},$
		 $ \frac{\theta^9+27\delta\theta^5+729\theta}{3^{7}},$ $\frac{\theta^{10}+27\delta\theta^6+729\theta^2}{3^8},$ $\frac{\theta^{11}+27\delta\theta^7+729\theta^3}{3^{9}}$\}, $\delta\in\{-1,1\}$\\
		%\cline{1-1}
		%\cline{3-5}
		B8 &$9$&$m_3\not\equiv\delta~(9)$& $48$&\tiny \{$1$, $\theta,$ $\frac{\theta^2}{3}, $ $\frac{\theta^3}{3^2},$ $\frac{\theta^4-54\delta}{3^3},$ $\frac{\theta^5-54\delta\theta}{3^4},$ $\frac{\theta^6-54\delta\theta^2}{3^4},$ $\frac{\theta^7-54\delta\theta^3}{3^5},$ $\frac{\theta^8+27\delta\theta^4+729}{3^6},$
		$ \frac{\theta^9+27\delta\theta^5+729\theta}{3^{7}},$ $\frac{\theta^{10}+27\delta\theta^6+729\theta^2}{3^8},$ $\frac{\theta^{11}+27\delta\theta^7+729\theta^3}{3^{8}}$\}, $\delta\in\{-1,1\}$\\
		\hline
	\end{longtable}
\end{center}
\end{theorem}\vspace{-0.3in}
%\begin{corollary}\label{cor1}
%	\indent Let $K=\Q(\theta)$ be an algebraic number field defined by a root $\theta$ of an irreducible polynomial $f(x)=x^{12}-2^{r}3^{s}$, where $r\in \{2,4,6,8,10\}$ and $s\in \{3,6,9\}$. Then 
%	\begin{equation*}
%		d_K=12\times2^{11(r+2)-2v_2(\ind\theta)}\times3^{11(s+1)-2v_3(\ind\theta)}
%	\end{equation*}
%\end{corollary}
\begin{theorem}\label{Th1.3}
	Let $K=\Q(\theta)$ be an algebraic number field with $\theta$ a root of an irreducible polynomial $f(x)=x^{12}-m$, where $m$ is a $12$-th power free integer. Suppose $d=\gcd(v_p(m),12)$ and $p$ is a prime number. Then the following hold:\begin{enumerate}
		\item  If   $p\mid m$ and $p\nmid v_p(m)$, then $$v_p(\ind\theta)=\frac{1}{2}[11(v_p(m)-1)+(d-1)]$$ and a $p$-integral basis of $K$ is 
 $$S=\left\{\frac{\theta^i}{p^{\lfloor{\frac{iv_p(m)}{12}\rfloor}}} |~ 1  \le i\le 12\right\}.$$  
\item If $p=3$ and $p\nmid m$, then  \[v_3(\ind\theta)=  \left\{ 
\begin{array}{cc}
	4, & \text{if}~ m^2\equiv1~(9),\\
	0, & \text{if}~ m^2\not\equiv 1~(9)\\
\end{array}
\right. \]
 and a $3$-integral basis of $K$ is
\[\left\{ 
\begin{array}{cc}
	1,\theta,\theta^2,\cdots,\theta^7, \frac{h(\theta)}{3}, \frac{\theta h(\theta)}{3}, \frac{\theta^2 h(\theta)}{3},\frac{\theta^3h(\theta)}{3}, & \text{if}~ m^2\equiv1~(9),\\
	1,\theta,\theta^2,\cdots,
\theta^{10},\theta^{11}, & \text{if}~ m^2\not\equiv 1~(9),\\
	
\end{array}
\right.\]
where $h(\theta)=\theta^8+m\theta^4+1.$

\item  If $p=2$ and $p\nmid m$, then 

\[v_2(\ind\theta)=  \left\{ 
\begin{array}{cc}
	
	9, & \text{if}~ m\equiv 1~(8),\\
	6, & \text{if}~ m\equiv 5~(8),\\
	0, & \text{if}~ m\equiv 3~(4)\\
	\end{array}
\right. \] and a $2$-integral basis of $K$ is 
\[\left\{ 
\begin{array}{cc}
	1, \theta, \theta^2, \theta^3, \theta^4, \theta^5, \frac{\theta^6-1}{2}, \frac{\theta(\theta^6-1)}{2}, \frac{\theta^2(\theta^6-1)}{2}, \frac{g(\theta)}{2^2}, \frac{\theta g(\theta)}{2^2}, \frac{\theta^2g(\theta)}{2^2}, & \text{if}~ m\equiv1~(8),\\
	1, \theta, \theta^2, \theta^3, \theta^4, \theta^5, \frac{\theta^6-1}{2}, \frac{\theta(\theta^6-1)}{2}, \frac{\theta^2(\theta^6-1)}{2}, \frac{\theta^3(\theta^6-1)}{2}, \frac{\theta^4(\theta^6-1)}{2}, \frac{\theta^5(\theta^6-1)}{2}, & \text{if}~ m\equiv5~(8),\\
	1,\theta,\theta^2,\cdots,\theta^{10},\theta^{11}, & \text{if}~ m\equiv 3 ~(4),\\
	
	\end{array}
\right.\] where $g(\theta)=\theta^9+\theta^6+\theta^3+1.$
\item If $p\nmid 12m$, then $v_p(\ind\theta)=0$ and $\{1,\theta, \theta^2,\cdots,\theta^{10}, \theta^{11}\}$ is a $p$-integral basis of $K.$
\end{enumerate}
  \end{theorem}

%\begin{theorem}\label{Th1.3}
%	Let $K=\Q(\theta)$ be same as in Theorem \ref{Th1.2}. Suppose  $p$ is a rational prime and $d=\gcd(v_p(m),12)$, then the following holds:\\
	%\vspace{-0.5cm}
	%\begin{enumerate}
	%\item If $p\mid m$ and $p\nmid d$, then $$v_p(\ind\theta)=\frac{1}{2}[11(v_p(m)-1)+(d-1)]$$.\\
	%\item  If $p\in\{2,3\}$ and $p\nmid m$, then 
	
		%\[v_2(\ind\theta)=  \left\{ 
		%\begin{array}{cc}
			
			%	9 & if~ m\equiv 1~(8)\\
		%6 & if~ m\equiv 5~(8)\\
		%0 & if~ m\equiv 3~(8)\\
		
		%\end{array}
	%\right. \] and 
%	\[v_3(\ind\theta)=  \left\{ 
%	\begin{array}{cc}
		%4 & if~ m^2\equiv1~(9)\\
		%0 & if~ m^2\not\equiv 1~(9)\\
	%\end{array}
	%\right. \] 
%\end{enumerate}
	
%\end{theorem}
The following corollary follows immediately from the above theorem.

\begin{corollary}\label{abc}
\indent Let $m$ be a square-free integer and $K=\Q(\theta)$ be an algebraic number field defined by a complex root $\theta$ of an irreducible polynomial $x^{12}-m$. Then $\{1,\theta,\theta^2,\cdots,\theta^{10},\theta^{11}\}$ is an integral basis of $K$ if and only if one of the following hold:
\vspace{-0.2cm}
\begin{enumerate}
	\item $2\nmid m$, $ 3\nmid m$, $m\equiv 3~(4)$ and $m^2\not\equiv 1~(9). $
	\item $2\nmid m$, $ 3\mid m$ and $m\equiv 3~(4)$.
	\item $2\mid m$, $ 3\nmid m$ and $m^2\not\equiv 1~(9).$
	\item $6\mid m$.
\end{enumerate} 

\end{corollary}

 \subsection{Construction of an explicit integral basis from $p$-integral basis}
 Let $\xi$ is a root of a monic irreducible polynomial over $\Q$ and  $L=\Q(\xi).$ It is well known that there exists an integral basis $\mathcal B:=\{\beta_0,\cdots ,\beta_{n-1}\}$ of $L$ such that $$\beta_0=1,~\beta_i=\frac{a_{i0}+a_{i1}\xi+\cdots +a_{i(i-1)}\xi^{i-1}+\xi^i}{d_i}$$ 
 with $a_{ij},~d_i$ in $\Z$ and the positive integer $d_i$ dividing $d_{i+1}$ for $1\leq i\leq n-1$; moreover  $[O_L : \Z[\xi]]=\displaystyle\prod_{i=1}^{n-1}d_i$  and the numbers $d_i$ are uniquely determined by $\xi$  ( \cite[Chapter 2, Theorem 13]{13}). Fix a prime $p$ and let $l_i=v_p(d_i)$, then $v_p(\ind\xi)=l_1+l_2+\cdots +l_{n-1}.$ As every integral basis of $L$ is its $p$-integral basis,  so is $\mathcal B_p^*:=\{1,\frac{\beta_1d_1}{p^{l_1}},\cdots ,\frac{\beta_{n-1}d_{n-1}}{p^{l_{n-1}}}\}$. It can be easily seen that if $\mathcal C:=\{1,\gamma_1,\cdots ,\gamma_{n-1}\}$ is another $p$-integral basis of $L$ where $\gamma_i$'$s$ are of the form  $$\gamma_i=\frac{c_{i0}+c_{i1}\xi+\cdots +c_{i(i-1)}\xi^{i-1}+\xi^i}{p^{k_i}} $$
 with $c_{ij},~k_i$ in $ \Z$ for $1\leq i\leq n-1$, then on writing each member of $\mathcal B_p^*$ as a $\Z_{(p)}$-linear combination of members of $\mathcal C$ and vice versa, we see that $l_i=k_i~ \forall~ i$ and hence 
 \vspace*{-0.05in}\begin{equation}\label{eqn p2}
 	v_p(\ind\xi)=l_1+l_2+\cdots +l_{n-1}=k_1+k_2+\cdots +k_{n-1}.
 	\vspace*{-0.05in}\end{equation}
 These $p$-integral bases of $L$ with $p$ running over all primes dividing $\ind\xi,$ quickly lead to a construction of an integral basis of $L$.\\
 \indent The following theorem proved in \cite[Theorem 2.2]{a}, describes a procedure for constructing an integral basis from all $p$-integral bases. We omit its proof.
 
 \begin{theorem}\label{Th 2.4} Let $L=\Q(\xi)$ be an algebraic number field of degree $n$ with $\xi$ an algebraic integer.  Let $\{\alpha_{r0},\alpha_{r1},\cdots,\alpha_{r(n-1)}\}$ be a $p_r$-integral basis of $L$, $1 \leq r\leq s $ with $\alpha_{r0}=1$, $\alpha_{ri}=\frac{c_{i0}^{(r)}+c_{i1}^{(r)}\xi+\cdots+c_{i(i-1)}^{(r)}\xi^{i-1}+\xi^i}{p_r^{k_{i,r}}}$, $1\leq i \leq n-1$, where $c_{ij}^{(r)}$ and $ 0\leq k_{i,r}\leq k_{i+1,r}$ are integers. If  $c_{ij} \in \Z$ are such that $c_{ij} \equiv c_{ij}^{(r)}($ $p_r^{k_{i,r}})$ for $1 \leq r\leq s $ and if $t_i$ stands for $\displaystyle \prod_{r=1}^{s} p_r^{k_{i,r}}$, then $\{\alpha_0,\alpha_1,\cdots,\alpha_{n-1}\}$ is an integral basis of $L$ where $\alpha_0=1$, $\alpha_i=\frac{c_{i0}+c_{i1}\xi+\cdots+c_{i(i-1)}\xi^{i-1}+\xi^i}{t_{i}}$ for $1\leq i\leq n-1$.
 \end{theorem}
 %\begin{proof} Clearly $\alpha_i$ is integral over $\Z_{(p)}$ for each prime $p$ not belonging to the set $\{p_1,\cdots ,p_{s}\}$ for  $1\leq i \leq n-1$. Keeping in mind  the choice of $c_{ij},$ we see that each $\alpha_i$ is integral over $\Z_{(p)}$ for $p$ belonging to the set $\{p_1,p_2,\cdots ,p_s\}$ and hence $\alpha_i$ is integral over $\Z$. Therefore if $\Gamma$ denotes the   subgroup of $\C$ defined by $\Gamma =\Z\alpha_0+\Z\alpha_1+\cdots +\Z\alpha_{n-1}$, then $\Z[\xi]\subseteq\Gamma \subseteq A_L$. Further by virtue of a basic result ( \cite[Chapter 2, Section 2, Theorem 2]{Bor} ),  the index of the subgroup $\Z[\xi]$ in $\Gamma$  equals $\displaystyle\prod_{i=1}^{n-1}t_i$. In view of (\ref{eqn p2}), we have$v_{p_{r}}(ind~\xi)=k_{1,r}+\cdots +k_{n-1,r} $ for $1\leq r\leq s$. So\centerline{$ind~\xi=\displaystyle\prod_{r=1}^{s}p_r^{k_{1,r}+\cdots +k_{n-1,r} }=\displaystyle\prod_{i=1}^{n-1}(\displaystyle\prod_{r=1}^{s}p_r^{k_{i,r}})= \displaystyle\prod_{i=1}^{n-1}t_i.$} Since $\Gamma \subseteq A_L$ and $[\Gamma:\Z[\xi]]=[A_L:\Z[\xi]]$, it follows that $A_L=\Gamma$ as asserted.\end{proof}
 
 We now provide some examples.
 \begin{example}\label{s}
 	Let $K=\Q(\theta)$ be an algebraic number field defined by a root $\theta$ of $f(x)=x^{12}-2352$. As $f(x)$ satisfies Eisenstein criterion with respect to $3$, so it is an irreducible polynomial. In view of Equation \eqref{eq2}, $D_f=-2^{68}3^{23}7^{22}.$ Applying Theorem \ref{Th1.2}, the set \{ $1, \theta, \theta^2, \frac{\theta^3}{2}, \frac{\theta^4}{2},
 	\frac{\theta^5}{2}, \frac{\theta^6+4\theta^3+12}{2^2}, \frac{\theta^7+4\theta^4+12\theta}{2^2}, \frac{\theta^8+4\theta^5+12\theta^2}{2^3},$ $\frac{\theta^9+2\theta^6+4\theta^3+8}{2^3}, \frac{\theta^{10}+2\theta^7+4\theta^4+8\theta}{2^4},$
 	$\frac{\theta^{11}+2\theta^8+4\theta^5+8\theta^2}{2^4}\}$ is a  $2$-integral basis of $K$ and therefore, we see that \{$1$, $\theta,$ $ \theta^2,$ $\frac{\theta^3}{2},$ $ \frac{\theta^4}{2},$
 	$\frac{\theta^5}{2},$ $\frac{\theta^6}{2^2},$ $\frac{\theta^7}{2^2},$ $\frac{\theta^8+4\theta^5+4\theta^2}{2^3},$ $\frac{\theta^9+2\theta^6+4\theta^3}{2^3},$ $\frac{\theta^{10}+2\theta^7+4\theta^4+8\theta}{2^4},$
 	$\frac{\theta^{11}+2\theta^8+4\theta^5+8\theta^2}{2^4}$\} is a $2$-integral basis of $K.$ Using Theorem \ref{Th1.3}, it follows that  
 	$\{1,\theta, \theta^2,\cdots,\theta^{10}, \theta^{11}\}$ is a $3$-integral basis of $K$ and  
 	$\{ 1,\theta,\theta^2, \theta^3,$ $\theta^4, \theta^5, \theta^6, \frac{\theta^7}{7}, \frac{\theta^8}{7},$
 	$ \frac{\theta^9}{7}, \frac{\theta^{10}}{7}, \frac{\theta^{11}}{7}\}$ is a $7$-integral basis of $K$. Let $p_1=2, p_2=3$ and $p_3=7.$ 
 According to the notations of  Theorem \ref{Th 2.4}, we see that $r=3,$ and we can take  $c_{82}=c_{85}=28$ (as $28\equiv 4 ~(8)$ and $28\equiv0~ (7)$), $c_{93}=28$. Also $42\equiv 2~(8)$, $56\equiv 8~(16)$, $-28\equiv 4~(16)$ and $-14\equiv2 ~(16)$ implies that we can choose  $c_{96}=42,$  $c_{10,1}=c_{11,2}=56, c_{10,4}=c_{11,5}=-28,$ $c_{10,7}=c_{11,8}=-14$, otherwise take $c_{ij}=0.$ Thus in view of \ref{Th 2.4}, we conclude that the  set 
 	\{ $1$,$\theta,$ $\theta^2,$  $\frac{\theta^3}{2},$ $ \frac{\theta^4}{2},$ $\frac{\theta^5}{2},$ $\frac{\theta^6}{4},$ $\frac{\theta^7}{28},$ $\frac{\theta^8+28\theta^5+28\theta^2}{56},$
 	$\frac{\theta^9+42\theta^6+28\theta^3}{56},$ $\frac{\theta^{10}-14\theta^7-28\theta^4+56\theta}{112},$   $\frac{\theta^{11}-14\theta^8-28\theta^5+56\theta^2}{112}$\} is an integral basis of $K.$
 \end{example} 
 
 \begin{example}\label{s}
 	Let $K=\Q(\theta)$ be an algebraic number field defined by $f(x)=x^{12}-60$. In view of Eisenstein criterion with respect to $3$, $f(x)$ is an irreducible polynomial. By  using Theorem \ref{Th1.2}, $v_2(\ind\theta)=15$ and a $2$-integral basis is given by 
 	$\{ 1, \theta, \theta^2, \theta^3, \theta^4, \theta^5, \frac{\theta^6-2\theta^3+2}{2^2},$ $ \frac{ \theta^7-2\theta^4+2\theta}{2^2},$
 	$\frac{\theta^8-2\theta^5+2\theta^2}{2^2}, \frac{\theta^9-2\theta^6+2\theta^3}{2^3},$$ 
 	$$\frac{\theta^{10}-2\theta^7+2\theta^4}{2^3}, \frac{\theta^{11}-2\theta^8+2\theta^5}{2^3}\}.$ By Theorem \ref{Th1.3}, we have $v_p(\ind\theta)$  $= 0$ and a $p$-integral basis is $\{1, \theta, \theta^2, \cdots \cdots  ,\theta^{10},  \theta^{11}\}$ for any odd prime $p.$ Thus using Equation \eqref{eq2} and the fact that $D_f=[O_K:\Z[\theta]]^2d_K,$ it follows that   $d_K=-2^{16}3^{23}5^{11}$.  Using Theorem \ref{Th 2.4}, we see that the set \{ $1, \theta, \theta^2, \theta^3, \theta^4, \theta^5, \frac{\theta^6-2\theta^3+2}{2^2}, \frac{ \theta^7-2\theta^4+2\theta}{2^2},
 	\frac{\theta^8-2\theta^5+2\theta^2}{2^2},$ $\frac{\theta^9-2\theta^6+2\theta^3}{2^3}, 
 	\frac{\theta^{10}-2\theta^7+2\theta^4}{2^3},$ $\frac{\theta^{11}-2\theta^8+2\theta^5}{2^3}\}$ is an integral basis of $K.$
 \end{example} 
 
 \begin{example}\label{s}
 	Let  $K=\Q(\theta)$ be an algebraic number field defined by a  root $\theta$ of an irreducible polynomial $f(x)=x^{12}-6a,$ where $a$ is a square free integer not divisible by $6.$ Then in view of Corollary \ref{abc}, $\{1,\theta,\theta^2,\cdots,\theta^{10}, \theta^{11}\}$ is an integral basis of $K.$ 
 \end{example}

 \section {Preliminary Results}
The following proposition to be used in the sequel follows immediately from Section 2.1.
\begin{prop}\label{prop1}
	Let $L=\Q(\eta)$ be an algebraic number field of degree $n$ with $\eta$ an algebraic integer and $p$ be a rational prime. Let $\beta_1,\beta_2,\cdots,\beta_{n-1}$ are $p$-integral elements of $L$ of the type $\beta_i=\frac{\eta^i+\displaystyle\sum_{j=1}^{i-1}c_{i,j}\eta^j}{p^{k_i}}$ where $c_{i,j}, k_i\in \Z$ with $0\le k_i \le k_{i+1}$ for $1\le i\le n-1.$ Then $\{1,\beta_1,\cdots,\beta_{n-1}\}$ is a $p$-integral basis of $L$ if and only if $v_p(\ind\eta)=\displaystyle\sum_{i=1}^{n-1}k_i,$ in which the integers $k_1,\cdots,k_{n-1}$ are uniquely determined by the prime $p$ and the element $\eta$ of $L.$ Moreover there always exists a $p$-integral basis of $L$ of the above type.
\end{prop}
 \indent \textbf{Gauss valuation, Newton polygons of first order and second order.} Throughout the paper, $\F_p$ denotes the finite field with $p$ elements and $Z_p$ denotes the ring of $p$-adic integers. Also $\bar{a}$ stands for the image of $a$ under the canonical homomorphism from $\Z_p$ onto $\F_p.$ 
  \begin{definition}\label{A} 
The Gauss valuation  of the field $\Q_p(x)$ of rational functions in an indeterminate $x$ which extends the valuation $v_p$ of $\Q_p$ and is defined on $\Q_p[x]$ by \begin{equation*}\label{Gau}
 v_{p,x}( a_0+a_1x+a_2x^2+.....+a_sx^s)= \displaystyle\min_{1\leq i\leq s} \{v_p(a_i)\}, ~a_i\in \Q_p.
  \end{equation*}
\end{definition}
 \begin{definition}\label{B}
 	Let $p$ be a prime number and $g(x) = x^n+a_{n-1}x^{n-1}+\cdots+a_0$ with $a_0\ne0$ be a polynomial
 	over $\Z_p$. To each non-zero term $a_ix^i$, we associate a point $(i, v_p(a_i))$ and form the set  $P =\{(i, v_p(a_i)) : 0\le i\le n, a_i\ne0 \}$. The  $p$-Newton polygon of $g(x)$ of first order (also called Newton polygon of $g(x)$ with	respect to $p$) is the polygonal path formed by the lower edges along the convex hull of points of $P$. Note that the slopes of the edges are increasing when calculated from left to right.
 	\end{definition}
 \begin{definition}\label{C}
 	 Let $g(x)=x^n+a_{n-1}x^{n-1}+\cdots+a_0$ be a polynomial over $\Z_p$ such that
 	the $p$-Newton polygon of $g(x)$ of first order consists of a single edge having negative slope, say $\lambda\in Q$. Let $\lambda=-\frac{h}{e}$, with $h$ and $e$ are coprime positive integers. Then we associate with $g(x)$ a polynomial $T_g(Y )\in \F_p[Y ]$ not
 	divisible by $Y$ of degree $\frac{n}{e}=d$ (say) defined by
 	\begin{equation}\label{3}
 		T_g(Y)=Y^d+\sum_{j=0}^{d-1}\overline{\Big(\frac{a_{ej}}{p^{v_p(a_0)+ej\lambda}}\Big)}Y^j
 	\end{equation}
 	The polynomial $T_g(Y )$ is called the residual polynomial of $g(x)$ with respect
 	to $p$.
 	\end{definition}
We now state the following weaker version of the theorem proved by Ore \cite{11} in a more general set up. Its proof is omitted.
\begin{theorem}\label{SK}
Let $p$ be a prime number. Let $L=\Q(\gamma)$ where $\gamma$ is a root of an irreducible polynomial $g(x)=x^n+a_{n-1}x^{n-1}+\cdots+a_0\in\Z[x]$, $a_0\ne0 $ with $g(x)\equiv x^n\mod p$. Suppose that
the $p$-Newton polygon of $g(x)$ of first order consists of a single edge with negative slope $\lambda$. If the residual polynomial $T_g(Y)\in \F_p[Y]$ of $g(x)$ associated to this edge is separable, then $v_p(\ind\gamma)$ equals the number of points with positive integer coordinates lying on or below the $p$-Newton polygon of $g(x)$.
\end{theorem}
\begin{definition}
	Let $p$,$L=\Q(\gamma)$,$g(x)$, $T_g(Y)$ be as in the above theorem. We say that $g(x)$ is $p$-regular, if the residual polynomial $T_g(Y)\in  \F_p[Y]$ of $g(x)$ has no repeated roots.
	\end{definition}
 Let $L=\Q(\gamma)$ where $\gamma$ is a root of a monic polynomial $g(x)=x^n+a_{n-1}x^{n-1}+\cdots+a_0\in\Z[x],~a_0\neq 0$. Let $p$ be a prime number such that  $g(x)\equiv x^n\mod p$. Suppose that the $p$-Newton polygon of $g(x)$ of first order consists of a single edge with slope $\lambda=\frac{-h}{e},$ where $h$ and $e$ are coprime positive integers and $e>1$. Suppose $T_g(Y)=\psi(Y)^s$ in $\F_p[Y]$, where $s\geq 2$ and $\psi(Y)$  is a monic irreducible polynomial over $\F_p$.
 In this case, we construct a key polynomial $\Phi(x)$ attached with the slope $\lambda$ satisfying the following conditions: 
 \begin{itemize}
\item[(i)] $\Phi(x)\equiv x^k\mod p$, for some natural number $k$.
\item[(ii)] The $p$-Newton polygon of $\Phi(x)$ of first order is one-sided with slope $\lambda$.
\item[(iii)] The residual polynomial of $\Phi(x)$  with respect to $p$ is $\psi(Y)\in\F_p[Y]$.
\item[(iv)] $\deg \Phi(x)= e\deg \psi(Y)$.
\end{itemize}
As mentioned in \cite[Section 2.2]{7}, the data $(x;~\lambda,~\psi(Y))$ defines a $p $-adic valuation $V$ on the field $\Q_p(x)$ with  $V(x)=h$, $V(p)=e$ and $V(\Phi(x))=he\deg(\psi(Y)$). If $p(x)=\displaystyle\sum_{0\leq i}^{}b_ix^i\in\Z_p[x]$ is any polynomial, then 
\vspace{-0.2in}\begin{equation}\label{eq3}
	V(p(x))=e \displaystyle\min_{0\leq i}^{}\{v_p(b_i)+i \lvert\lambda\rvert\}.
\end{equation}
We define the above valuation $V$ to  the valuation of second order.
If $g(x)=\displaystyle\sum_{i=0}^{u} a_i(x)\Phi(x)^i$ is a $\Phi$-adic expansion of $g(x)$ in $\Z_p[x]$, then the $V$-Newton polygon of $g(x) $ of second order (also called  $V$-Newton polygon of $g(x)$) is the lower convex hull  of the points of the set $\{(i, V(a_{i}(x)\Phi(x)^{i})), 0\le i\le u$\} of the  Euclidean plane.\\
\indent Let the $V$-Newton polygon of $g(x)$ of second order has $k$-edges, say $E_1,\cdots , E_k$, with negative slopes $\lambda_1,\cdots,\lambda_k$.  Let $\lambda_t=\frac{-h_t}{e_t}$, where 
$h_t$ and $e_t$ are coprime positive integers and $l_t$ denote the projection
to the horizontal axis of the side of slope $\lambda_t$ for $1\leq t \leq k.$ Then, there is a natural residual polynomial $\psi_t(Y)$ of second order attached to each edge $E_t,$ whose degree coincides
with the degree of the edge (i.e. $\frac{l_t}{e_t}$) \cite[Section 2.5]{7}. Only those integral points of the
$ V $-Newton polygon of $g(x)$ which lie on the edge, determine a non-zero coefficient of this
second order residual polynomial.
We define $g(x)$ to be $\psi_t $-regular when the second order residual polynomial $\psi_t(Y)$ attached to the side $E_t$ of the $V$-Newton polygon of $g(x)$ of second order is separable in
$\frac{\F_p[Y ]}{\langle\psi_t(Y)\rangle}$. We define $g(x)$ to be $V$-regular if $g(x)$ is $\psi_t$-regular for each $t$, $1\leq t\leq k$.
\begin{definition}\label{D2.6}
Let $\Phi(x)\in \Z[x]$ be a monic polynomial, and let
$f(x)=a_n(x)\Phi(x)^n +\cdots+ a_1(x)\Phi(x)+a_0(x)$,
with $a_i(x)\in \Z[x]$, $\deg a_i(x)<\deg \Phi(x)$, be the $\Phi$-adic expansion of $f(x)$. Then we define
the quotients attached to this $\Phi$-expansion, by definition, the different quotients
$q_1(x),\cdots, q_n(x) $ that are obtained along the computation of the coefficients of the expansion:\vspace{-0.1in}
$$f(x)=\Phi (x)q_1(x)+a_0(x),$$
$$q_i(x)=\Phi (x)q_{i+1}(x)+a_i(x),~\forall 1\le i\le n.$$
\end{definition}
The following theorem is weaker version of the theorems proved by  Guardia, Montes and Nart in $2012$ [\cite{7}, Theorem $4.18$] and in $2015$ \cite{8}. 
\begin{theorem}\label{ABC}
Let $p$ be a prime number. Let $L=\Q(\gamma)$ where $\gamma$ is a root of a monic polynomial $g(x)=x^n+a_{n-1}x^{n-1}+\cdots+a_0\in \Z[x], a_0\ne 0$ with $g(x)\equiv x^n\mod p$. Suppose that the $p$-Newton polygon of $g(x)$ of first order consists of a single edge of negative slope $\lambda=-\frac{h}{e}$ with $\gcd(h,e)=1$ and the residual polynomial of $g(x)$ is given by $T_g(Y)=\psi_1(Y)^{r_1}\psi_2(Y)^{r_2}\cdots\psi_s(Y)^{r_s},$  is product of powers of irreducible polynomials $\psi_i(Y)\in \F_p[Y]$, where $\psi_i(Y)\ne Y,$ $r_i\ge 2$ for $1\le i \le s$. Let $\Phi_i(x)$ be the key polynomial attached to $\lambda$ and $V_i$ be the corresponding second order valuation determined by $(x;\lambda,\psi_i(Y))$ for all $i =1,2,\cdots ,s$. For every $1\le i\le s$, let  $\mu_i=\deg\psi_i(Y)$ and $V_i$-Newton polygon of $g(x)$ has $z_i$ edges $E_{i1},\cdots ,E_{iz_i}$ of  negative slopes
$\lambda_{i1},\cdots,\lambda_{iz_i}$. If $g(x)$ is $V_i$ -regular for each $i=1,2,\cdots,s$, then the following hold:\vspace{-0.1in}
\begin{enumerate}\item $v_p(\ind\gamma)= N_1+\displaystyle\sum_{i=1}^{s}\mu_iN^i_2$, where
$N_1$ is the number of points with positive integer coordinates lying on or below the $p$-Newton polygon of $g(x)$ and $N^i_2$ denote the number of points with positive integer coordinates lying on or below the $V_i$ -Newton polygon of $g(x)$ and lying above the horizontal line passing through the last vertex of this polygon.
\item Let $y_j$ denote the ordinate of the point of the $p$-Newton polygon of $g(x)$ of first order
with abscissa $j$. Then the  set $S=\displaystyle\cup_{i=1}^sS_i$ is a $p$-integral basis
of $L,$ where for each $i=1,2,\cdots,s,$
\begin{equation}
	 S_i =\left
\{\frac{\theta^{n-u}q_{ij}(\theta)}{p^{\lfloor{y_u+\frac{Y_{ij}-jV(\Phi_i(x))}{e}}\rfloor}}~:~ n-e\mu_i < u\le n, b_{it}-e_{it} f_{it} < j\le b_{it}, 1\le t\le z_i\right \},
\end{equation} with $q_{ij}(\theta)$ is the $j$-th quotient
in the $\Phi_i$-adic expansion of $ g(x)$ as in Definition \ref{D2.6}, $\lambda_{it}=\frac{-h_{it}}{e_{it}}$
with $\gcd(h_{it}, e_{it})=1,$ $[a_{it},b_{it}]$ denote the projection to the horizontal axis of the edge  $E_{it},$ $f_{it} =\frac{ b_{it}-a_{it}}{e_{it}}$
for $1\le t \le z_i$ and  $Y_{ij}$ denote the ordinate of the point of the $V_i$-Newton polygon
of $g(x)$ of second order with abscissa $j.$
\end{enumerate}

\end{theorem}

The elementary lemma stated below is well known (see \cite[Problem 435]{10}).
\begin{lemma}\label{lemma1}
	Let $t$, $b$ be positive integers with $\gcd(t, b)=c$. Let $P$ denote the set of points in the plane with positive integer coordinates lying inside or on the triangle with vertices  $(0,0), (t, 0), (0, b)$. Then 
		$\#P=\sum_{i=0}^{t-1}\lfloor\frac{ib}{t}\rfloor=\frac{1}{2}[(t-1)(b-1)+(c-1)],$ where $\#P$ is cardinality of the set $P$. 
\end{lemma}

 \section{ Proof of Theorems \ref{Th1.2}}
  
\begin{proof}[\textbf{Proof of Theorem \ref{Th1.2}}]
 
In what follows, $N_1$ will stand for the number of points with positive integer coordinates lying on or below the $p$-Newton polygon of $f(x)$ of first order
 and $N_2$ stands for  the number of points with positive integer coordinates lying on or below the $V$-Newton polygon of $f(x)$ of second order and lying above the horizontal line passing through the last vertex of the polygon.\\ 
\textbf{Case A1:} $v_2(m)=2$, $m_2\equiv 1~(4)$. In this case, $f(x)\equiv x^{12}~(2)$. The $2$-Newton polygon of $f(x)$ of first order has a single edge, say $S$, joining the points $(0,
~2)$ and $(12,~0)$ having slope $\lambda=-\frac{1}{6}$. The residual polynomial of $f(x)$  associated to $S$ is $Y^2+\bar{1}=(Y+\bar{1})^2\in \F_2[Y]$, which is not separable. Therefore $f(x)$ is not $2$-regular. In view of Lemma \ref{lemma1}, we have  $N_1=6$. Set $\psi(Y)=Y+\bar{1}$.\\ 
\indent For the  second order Newton polygon, we  define some numerical invariants, $h=1$, $e=6$ and $\mu=\deg(\psi(Y))=1$, where $h$ and $e$ respectively are the numerator and denominator of $\lambda$.
Choose  $\Phi(x)=x^6+2$, then one can easily check that $\Phi(x)$ is a key polynomial attached with slope $\lambda$. In view of Equation  \eqref{eq3}, we  define the valuation $V$ of second order on $\Q_2(x)$ attached to the data $(x;\lambda,\psi(Y))$ such that $V(x)=h=1$, $V(2)=e=6$ and $V(\Phi)=he\mu=6$. The $\Phi$-expansion of $f(x)$ is  given by
 $f(x)=\Phi^2(x)-4\Phi(x)+4-m.$
\indent The  $V$-Newton polygon of $f(x)$ of second order is the lower convex hull of the points $(0,~6v_2(4-m))$, $(1,~18)$ and $(2,~12)$. If $m_2\equiv1~(8)$, then the  $V$-Newton polygon of $f(x)$  has two edges of negative slope (see Figure \ref{fig1}). The first edge, say $S_1$, is the line segment joining the point $(0,~6v_2(4-m)$ with $(1,~18)$ and the second edge, say $S_2$, is the one joining the point 
$(1,~18)$ to $(2,~12)$. For each $i=1,2$, the residual polynomial  associated with $S_i$ is linear. So $f(x)$ is $V$-regular. According to the notations of Theorem \ref{ABC}, we have $y_u=-\frac{1}{6}u+2$, where $7\le u\le 12$, $a_1=0$, $b_1=1$, $a_2=1$, $b_2=2$,  $e_1=e_2=1$, $f_1=f_2=1$, $Y_1=18$, $Y_2=12$, $q_1(\theta)=\theta^6-2 $ and $q_2(\theta)=1$.
\vspace{-0.06in}
\begin{figure}[H]
	\centering
	\begin{minipage}{0.45\textwidth}
		\begin{tikzpicture}[scale=0.7]
			\draw[thick, ->] (-0.5, 0) -- (7, 0);
			\draw[thick, ->] (0, -0.5) -- (0, 7.5);
			\draw[thick] (0, 7) .. controls (2,4) and  (2,4) ..(2,4);
			\draw[thick] (2, 4) .. controls (5,2) and  (5,2) ..(5,2);
			
			\draw[dashed] (0,2) .. controls (5,2) and  (5,2) ..(5,2);

			\draw (0, 7) node[right]{$(0,6v_2(4-m))$};
			\draw (-.35, 6.9) node[right]{$\bullet$};
			\draw (3.8,4.6) node[left]{$(1,18)$};
			\draw (5, 2) node[right]{$(2,12)$};
			
			\draw (0, 2) node[left]{$(0,12)$};
			\draw (5.3, 2) node[left]{$\bullet$};
			\draw (1.6, 2.17) node[right]{$\star$};
			\draw (1.6, 2.536) node[right]{$\star$};
			\draw (1.6, 2.902) node[right]{$\star$};
			\draw (1.6, 3.268) node[right]{$\star$};
			\draw (1.6, 3.634) node[right]{$\star$};
			\draw (1.6, 4.01) node[right]{$\star$};
			
		\end{tikzpicture}
		\caption{$m_2\equiv 1~(8);$ $V$-Newton polygon of $f(x)$.}
		\label{fig1}
	\end{minipage}\hfill
	\begin{minipage}{0.45\textwidth}
		\begin{tikzpicture}[scale=0.7]
			\draw[thick, ->] (-0.5, 0) -- (8, 0);
			\draw[thick, ->] (0, -0.5) -- (0, 7.5);
			\draw[thick] (0, 6.5) .. controls (7,2) and  (7,2) ..(7,2);
			
			\draw[dashed] (0,2) .. controls (7,2) and  (7,2) ..(7,2);

			\draw (0, 6.5) node[right]{$(0,24)$};
			\draw (2.7,4.6) node[left]{$(1,18)$};
			\draw (6, 1.3) node[right]{$(2,12)$};
			
			\draw (0, 2) node[left]{$(0,12)$};
			\draw (.35, 6.45) node[left]{$\bullet$};
			\draw (7.32, 2) node[left]{$\bullet$};
			\draw (2.7, 2.3) node[right]{$\star$};
			\draw (2.7, 2.74) node[right]{$\star$};
			\draw (2.7, 3.18) node[right]{$\star$};
			\draw (2.7, 3.62) node[right]{$\star$};
			\draw (2.7, 4.06) node[right]{$\star$};
			\draw (2.7, 4.52) node[right]{$\star$};
			
		\end{tikzpicture}
		\caption{$m_2\equiv 5~(8);$ $V$-Newton polygon of $f(x)$.}
		\label{fig1.1}
	\end{minipage}
	\end{figure}\vspace{-0.2in}
\noindent If $m_2\equiv 5~(8)$, then the $V$-Newton polygon of $f(x)$ has a single edge, say $S'$, of negative slope (see Figure \ref{fig1.1}). The edge $S'$ is the line segment joining the points $(0,~24)$ and $(2,~12)$ with a lattice point $(1,~18)$ lying on it. The residual polynomial  associated to $S'$  is $Y^2+Y+\bar{1}\in\F_2[Y]$. Thus $f(x)$ is $V$-regular.  
Here we see that $y_u=-\frac{1}{6}u+2,$ where $7\le u\le 12,$  $[a_1,b_1]=[0,2],$  $e_1=1,$  $f_1=2,$  $ Y_1=18,$  $Y_2=12,$  $q_1(\theta)=\theta^6-2 $ and $q_2(\theta)=1$. From Figures \ref{fig1} and \ref{fig1.1}, we get $N_2=6$. By Theorem \ref{ABC}, it follows that $v_2(\ind\theta)=N_1+\mu N_2=6+6=12$ and  the set 
$$ A_1=\left\{1, \theta, \theta^2, \theta^3, \theta^4, \theta^5, \frac{\theta^6-2}{2^2}, \frac{\theta^7-2\theta}{2^2},\frac{\theta^8-2\theta^2}{2^2}, \frac{\theta^9-2\theta^3}{2^2}, \frac{\theta^{10}-2\theta^4}{2^2}, \frac{\theta^{11}-2\theta^5}{2^2} \right\},$$ is a $2$-integral basis of $K$.\\\\
\indent Note that in the above case, if $m_2\equiv3\mod8$, then the $V$-Newton polygon of $f(x)$ of second order being the lower convex hull of the points $(0,~18)$, $(1,~18)$, $(2,~12)$ has a single edge of negative slope. The residual polynomial  associated to this edge is $Y^2+\bar{1}\in\F_2[Y],$ which is not separable. Therefore  we take a different key polynomial for this situation.\\\\
\noindent\textbf{Case A2:} $v_2(m)=2$, $m_2\equiv3~(8)$. Arguing as in Case A1, we observe that $f(x)$ is not $2$-regular. In this case  $\lambda$, $h$, $e$, $\mu$, $N_1$ and $\psi(Y)$ will be same as in the Case A1. Take $\Phi(x)=x^6+2x^3-2$. Clearly  $\Phi(x)$ is a key polynomial attached to $\lambda$. The data $(x;\lambda,\psi(Y))$ determines a $2$-adic valuation $V$ on $\Q_2(x)$ given in Equation \eqref{eq3} such that $V(x)=1$, $V(2)=6$ and $V(\Phi)=6$.
The $V$-Newton polygon of $$f(x)=\Phi^2(x)+(8-4x^3)\Phi(x)-16x^3+12-m$$ being the lower convex hull of the points $(0,~27)$, $(1,~21)$ and $(2,~12)$ has a single edge, say $S'$, of negative slope. The edge $S'$ is the line segment joining the points $(0,~27)$, and $(2,~12)$. The residual polynomial  associated to $S'$ is linear. Thus $f(x)$ is $V$-regular and $N_2=7$. By virtue of the Theorem \ref{ABC}, we have $v_2(\ind\theta)=N_1+\mu N_2=6+7=13.$
According to the notations of Theorem \ref{ABC},   $y_u=-\frac{1}{6}u+2$, where $ 7\le u\le 12,$  $Y_j=-\frac{15}{2}j+27$, for $j=1,2$, $[a_1,b_1]=[0,2]$,  $e_1=2,$ $f_1=1,$  $q_1(\theta)=\theta^6-2\theta^3+6$ and $q_2(\theta)=1$. Hence the set 
$$
 A_2=\left \{1, \theta, \theta^2, \theta^3, \theta^4, \theta^5, \frac{q_1(\theta)}{2^2}, \frac{\theta q_1(\theta)}{2^2},\frac{\theta^2 q_1(\theta)}{2^2},
 \frac{\theta^3 q_1(\theta)}{2^2},
 \frac{\theta^{4}q_1(\theta)}{2^2}, \frac{\theta^{5}q_1(\theta)}{2^3}\right\},$$ is a $2$-integral basis of $K$\\\\
 \indent One can easily  verify that when  $m_2\equiv7\mod8$, then the $V$-Newton polygon of $f(x)$ of second order being the lower convex hull of the points $(0,~24)$, $(1,~21)$, $(2,~12)$ has a single edge of negative slope. The residual polynomial associated to this edge has a repeated root. Therefore we choose some  different key polynomial.\\\\
\noindent\textbf{Case A3:} $v_2(m)=2$, $m_2\equiv7~(8)$. In this case  $\lambda$, $e$, $\psi(Y)$, $\mu$ and $N_1$ will be same as in Case A1. Proceeding same as in Case A1, we observe that $f(x)$ is not $2$-regular. Let $\Phi(x)=x^6+2x^3+2$, then   $V(\Phi)=6$.
The $V$-Newton polygon of $f(x)=\Phi^2(x)-4x^3\Phi(x)-4-m,$  is the lower convex hull of the points $(0,~6v_2(4+m))$, $(1,~21)$  and $(2,~12)$. If $m_2\equiv 7~(16)$, then it has a single  edge, say $S'$, joining the points $(0,~30)$ and $(2,~12)$ with a lattice point $(1,~21)$ lying on it. The residual polynomial  associated with $S'$ is $Y^2+Y+\bar{1}\in \F_2[Y]$. %Here we have $y_u=-\frac{1}{6}u+2,$ where $7\le u\le 12,$ $[a_1,b_1]=[0,2],$  $e_1=1,$ $f_1=2,$ $Y_1=21,$  $Y_2=12,$ $q_1(\theta)=\theta^6-2\theta^3+2$ and $q_2(\theta)=1$.
If $m_2\equiv15~(16)$, then the $V$-Newton polygon has two edges of negative slope. The first edge, say $S_1$, is the line segment joining the point $(0,~6v_2(m+4))$ with $(1,~21)$ and the second edge, say $S_2$, is the one joining the point $(1,~21)$ to $(2,~12)$. For each $i=1,2$, the residual polynomial of second order associated with $S_i$ is linear. %Here $y_u=-\frac{1}{6}u+2,$ where  $7<u\le 12,$ $[a_1,b_1]=[0,1],$ and $[a_2,b_2]=[1,2],$  $e_1=e_2=1,$ $f_1=f_2=1,$ $Y_1=21$ and $Y_2=12$.
 Thus Lemma \ref{lemma1} and Theorem \ref{ABC} implies that  $v_2(\ind\theta)=N_1+ \mu N_2=6+9=15$ and the  following set is a $2$-integral basis of $K.$
$$
 A_3=\left \{1, \theta, \theta^2, \theta^3, \theta^4, \theta^5, \frac{q_1(\theta)}{2^2}, \frac{\theta q_1(\theta)}{2^2},\frac{\theta^2 q_1(\theta)}{2^2},
 \frac{\theta^3 q_1(\theta)}{2^3},
 \frac{\theta^{4}q_1(\theta)}{2^3}, \frac{\theta^{5}q_1(\theta)}{2^3}\right\},$$ where $q_1(\theta)=\theta^6-2\theta^3+2.$\\\\
\textbf{Case A4:}  $v_2(m)=4$, $m_2\equiv3~(4)$. In this case $f(x)\equiv x^{12}~(2)$. The $2$-Newton polygon of $f(x)$ of first order has a single edge  joining the points  $(0,~4), $  and $(12,~0)$ with slope $\lambda=-\frac{1}{3}(=-\frac{h}{e})$. The residual polynomial associated to this edge  is $Y^4+\bar{1}=(Y+\bar{1})^4$, which has a repeated root. Therefore $f(x)$ is not $2$-regular. Let $\psi(Y)=Y+\bar{1}$. In view of Lemma \ref{lemma1}, we have $N_1=18$. Let $\Phi(x)=x^3-2$ and  $V\big(\displaystyle\sum_{i\geq 0}a_ix^i\big)=3 \min_{i\ge0}\{v_2(a_i)+i\frac{1}{3}\}.$ Then $V(x)=1$, $V(2)=3$ and $V(\Phi)=3$. The $V$-Newton polygon of
\begin{equation}\label{eq3.4}
	f(x)=(\Phi(x)+2)^4=\Phi^4(x)+8\Phi^3(x)+24\Phi^2(x)+32\Phi(x)+16-m
\end{equation}
 being the lower convex hull of the points  $(0,~15)$, $(1,~18)$, $(2,~15)$, $(3,~18)$ and $(4,~12)$ has a single  edge, say $S'$, of negative slope. The edge $S'$ is the line segment joining the points $(0,~15)$ and $(4,~12)$. The  residual polynomial  associated with $S'$ is linear. Thus $f(x)$ is $V$-regular. Clearly  $y_u=-\frac{1}{3}u+4$, where $10\le u\le 12$, $e=3,$ $a_1=0,$ $b_1=4,$ $e_1=4,$ $f_1=1,$ $ Y_j=-\frac{3}{4}j+15$, for $j=1,2,3,4$, $q_1(\theta)=\theta^9+2\theta^6+4\theta^3+8$, $q_2(\theta)=\theta^6+4\theta^3+12$ and $q_3(\theta)=\theta^3+6$. Thus by Theorem \ref{ABC}, we see that  $v_2(\ind\theta)=N_1+\mu N_2=18+3=21.$ and  
$1, \theta, \theta^2, 
\frac{q_3(\theta)}{2},\frac{\theta q_3(\theta)}{2},
\frac{\theta^2 q_3(\theta)}{2}, \frac{q_2(\theta)}{2^2},\frac{\theta q_2(\theta)}{2^2},\frac{\theta^2 q_2(\theta)}{2^3}, \frac{q_1(\theta)}{2^3}, \frac{\theta q_1(\theta)}{2^4}, \frac{\theta^2q_1(\theta)}{2^4}$ are $2$-integral elements of $K.$ Keeping in mind $\frac{q_3(\theta)}{2}=\frac{\theta^3}{2}+3,$ we see that $\frac{\theta^3}{2},\frac{\theta^4}{2}$ and $\frac{\theta^5}{2}$ are algebraic integers. Hence in view of Proposition \ref{prop1}, the  set 
$$A_4=\left\{1,\theta,\theta^2,\frac{\theta^3}{2},\frac{\theta^4}{2},\frac{\theta^5}{2},\frac{q_2(\theta)}{2^2},\frac{\theta q_2(\theta)}{2^2},\frac{\theta^2 q_2(\theta)}{2^3},\frac{q_1(\theta)}{2^3}, \frac{\theta q_1(\theta)}{2^4}, \frac{\theta^2q_1(\theta)}{2^4}\right\}$$  is a $2$-integral basis of $K.$\\
\textbf{Case A5:}  $v_2(m)=4$, $m_2\equiv5~(8)$. Here $f(x)\equiv x^{12}~(2)$. Arguing as in Case A4, we see that $f(x)$ is not $2$-regular and $N_1=18$. Let $\Phi(x)$ and $V$ be as in the previous case. Keeping in mind the $\Phi$-expansion of $f(x)$ given in \eqref{eq3.4}, one can check that the  $V$-Newton polygon of $f(x)$ has a single edge, say $S'$, of negative slope. The edge $S'$ is the line segment joining the points $(0,~18)$ and $(4,~12)$ with a lattice point $(2,~15)$ lying on it. The  residual polynomial  associated with $S'$ is $Y^2+Y+\bar{1}\in \F_2[Y]$, which has no repeated roots. Therefore $f(x)$ is $V$-regular. In this case $N_2=8$. Hence by Theorem \ref{ABC}, $v_2(\ind\theta)=N_1+ \mu N_2=26$ and the set 
$$
A_5=\left \{1, \theta, \theta^2, 
\frac{q_3(\theta)}{2},\frac{\theta q_3(\theta)}{2},
\frac{\theta^2 q_3(\theta)}{2^2},   
\frac{q_2(\theta)}{2^3}, \frac{\theta q_2(\theta)}{2^3}, \frac{\theta^2q_2(\theta)}{2^3}, \frac{q_1(\theta)}{2^4}, \frac{\theta q_1(\theta)}{2^4}, \frac{\theta^2q_1(\theta)}{2^5}\right\},$$ where 
$q_1(\theta)=\theta^9+2\theta^6+4\theta^3+8$, $q_2(\theta)=\theta^6+4\theta^3+12$ and $q_3(\theta)=\theta^3+6,$ is a $2$-integral basis of $K$.\\\\
\textbf{Case A6:}  $v_2(m)=4$, $m_2\equiv1~(8)$. Here $f(x)\equiv x^{12}~(2)$. Proceeding same as in Case A4, we observe that $f(x)$ is not $2$-regular and $N_1$=18. In this case $\lambda$, $\psi(Y)$, $\Phi(x)$ and $V$  are same as in Case A4. Keeping in mind \eqref{eq3.4}, the $V$-Newton polygon of $f(x)$  is the  lower convex hull of the points $(0,~3v_2(16-m))$, $(1,~18)$, $(2,~15)$, $(3,~18)$ and $(4,~12)$.   If $m_2\equiv1~(16)$, then it has three edges of negative slope. The first edge joins  $(0,~3v_2(16-m))$ with $(1,~18)$, the second edge joins $(1,~18)$ with $(2,~15)$  and the third edge is from $(2,~15)$ to $(4,~12)$. The  residual polynomial  attached to each side is linear.  %$y_u=-\frac{1}{3}u+4$, where $10\le u\le 12$, $e=3, [a_1,b_1]=[0,1], [a_2,b_2]=[1,2], [a_3,b_3]=[2,4], e_i=1$, for $1\le i\le2$, $e_3=2$, $f_i=1$, for $1\le i\le3$, $Y_1=18, Y_2=15, Y_3=\frac{27}{2}$ and  $Y_4=12$.
  If $m_2\equiv9~(16)$, then the  $V$-Newton polygon has two edges of negative slope. The  first edge, say $S_1$, is the line segment joining the points $(0,~21)$ and $(2,~15)$ with a lattice point $(1,~18)$ lying on it and the second edge, say $S_2$, joins $(2,~15)$ to $(4,~12)$. The residual polynomial  associated  to  $S_1$ and $S_2$ is $Y^2+Y+\bar{1}$ and $Y+\bar{1}$ respectively. Therefore $f(x)$ is $V$-regular. It is easy to check that $N_2=10$. % According to the notations of Theorem \ref{ABC}, we have  $y_u=-\frac{1}{3}u+4$, where $10\le u\le 12$, $e=3, a_1=0, b_1=2, a_2=2, b_2=4, e_1=1$, $e_2=2$, $f_1=2, f_2=1$, $Y_1=18, Y_2=15, Y_3=\frac{27}{2},$  $Y_4=12$ and  $v_2(\ind\theta)=N_1+\deg\psi(Y) N_2=28$.
   Hence in view of Theorem \ref{ABC}, $v_2(\ind\theta)=N_1+\deg\psi(Y) N_2=28$ and the following  set is a $2$-integral basis of $K$.
$$
A_6=\left \{1, \theta, \theta^2,
\frac{q_3(\theta)}{2},\frac{\theta q_3(\theta)}{2},
\frac{\theta^2 q_2(\theta)}{2^2}, \frac{q_2(\theta)}{2^3},\frac{\theta q_2(\theta)}{2^3},
\frac{\theta^2 q_2(\theta)}{2^3}, \frac{q_1(\theta)}{2^5}, \frac{\theta q_1(\theta)}{2^5}, \frac{\theta^2q_1(\theta)}{2^5} \right\},$$
where $q_1(\theta)=\theta^9+2\theta^6+4\theta^3+8$, $q_2(\theta)=\theta^6+4\theta^3+12$ and  $q_3(\theta)=\theta^3+6$.\\\\
\textbf{Case A7:} $v_2(m)=6$, $m_2\equiv1~(4)$. In this case $f(x)\equiv x^{12}~(2)$. The 2-Newton polygon of $f(x)$ of  has a single edge  joining  the points $(0,~6)$ and $(12,~0)$ with slope $\lambda=-\frac{1}{2}(=-\frac{h}{e})$. In view of Lemma \ref{lemma1}, we have $N_1=30$. The residual polynomial of $f(x)$  associated with the edge is $Y^6+\bar{1}=(Y+\bar{1})^2(Y^2+Y+\bar{1})^2\in\F_2[Y],$ which is not separable. Let $\psi_1(Y)=Y+\bar{1}$ and $\psi_2(Y)=Y^2+Y+\bar{1},$ then we have two types $(x;-\frac{1}{2},Y+\bar{1})$ and $(x;-\frac{1}{2},Y^2+Y+\bar{1})$. Now we examine each type separately.\\
\indent Consider the first type  $(x;\lambda,\psi_1(Y))$. Let   $\Phi(x)=x^2-2$ and the second order valuation $V_1$ on $\Q_2(x)$ attached to the data $(x;\lambda,\psi_1(Y))$ defined as $V_1\big(\displaystyle\sum_{i\geq 0}^{}a_ix^i\big)=2\min_{i\ge0}\{v_2(a_i)+i\frac{1}{2}\}$, then $V_1(x)=1$, $V_1(2)=2$ and $V_1(\Phi)=2$. In this case the $V_1$-Newton polygon of  
\begin{equation}\label{eq3.5}
	f(x)=\Phi^6(x)+12\Phi^5(x)+60\Phi^4(x)+160\Phi^3(x)+240\Phi^2(x)+192\Phi(x)+64-m
\end{equation}
  is the lower convex hull of the points $(0,~2v_2(64-m)$, $(1,~14)$, $(2,~12)$, $(3,~16)$, $(4,~12)$, $(5,~14)$ and $(6,~12)$. If $m_2\equiv1~(8)$, then it has two edges of negative slope. The first edge is the line segment joining the point $(0,2v_ 2(64-m))$ with $(1,14)$  and the second is the one joining  $(1,14)$ to $(2,~12)$. The residual polynomial   associated with each edge is linear. If $m_2\equiv5~(8)$, then the $V_1$-Newton polygon  has a single edge of negative slope joining the points $(0,~16)$ and $(2,~12)$ with a lattice point $(1,~14)$ lying on it. The  residual polynomial  associated to this edge is $Y^2+Y+\bar{1}\in \F_2[Y]$. Therefore $f(x)$ is $V_1$-regular and   $N^1_2=2$.\\
\indent Consider the second type  $(x;-\frac{1}{2},\psi_2(Y))$. Let $\Phi(x)=x^4+2x^2+4$ and  the second order valuation $V_2$ attached to $(x;\lambda,\psi_2(Y))$ on $\Q_2(x)$ defined as  $V_2\big(\displaystyle\sum_{i\geq 0}^{}a_ix^i\big)=2 \min_{i\ge0}\{v_2(a_i)+i\frac{1}{2}\},$ then $V_2(x)=1$, $V_2(2)=2$ and  $V_2(\Phi)=4$. The $V_2$-Newton polygon  of
\begin{equation}
	f(x)= \Phi^3(x)-6x^2\Phi^2(x)+(16x^2-32)\Phi(x)+64-m
\end{equation}
 is the lower convex hull of the points $(0,~2v_2(64-m))$, $(1,~14)$, $(2,~12)$ and $(3,~12)$.\\
If  $m_2\equiv 1~(8)$, then it has two edges of negative slope. The first edge is the line segment joining the  point $(0,~2v_2(64-m))$ with $(1,~14)$ and the second edge is the one joining the point $(1,~14)$ to $(2,~12)$. The residual polynomial attached to each edge is linear.
If $m_2\equiv5~(8)$, then $V_2$-Newton  polygon  has a single edge   joining the points $(0,~16)$ and $(2,~12)$ with a lattice point $(1,~14)$ lying on it. The  residual polynomial  associated to this edge  is a separable polynomial of degree $2.$ In this case  $N^2_2=2$.\\ 
\indent Keeping in mind Theorem \ref{ABC}, we see that $v_2(\ind\theta)=N_1+\deg \psi_1(Y)N_2^1+\deg\psi_2(Y)N_2^2=30+2+4=36$ and a $2$-integral basis of $K$ is 
{\footnotesize$$\left \{ 
\frac{q'_2(\theta)}{2^2},\frac{\theta q'_2(\theta)}{2^2},
\frac{\theta^2 q'_2(\theta)}{2^3}, \frac{\theta^3 q'_2(\theta)}{2^3}, \frac{ q_2(\theta)}{2^4}, \frac{\theta q_2(\theta)}{2^4}, \frac{q'_1(\theta)}{2^5},\frac{\theta q'_1(\theta)}{2^5},    
\frac{q_1(\theta)}{2^6}, \frac{\theta q_1(\theta)}{2^6}, \frac{\theta^2 q'_1(\theta)}{2^6}, \frac{\theta^3 q'_1(\theta)}{2^6} \right\},$$}
where 
$q_1(\theta)=\theta^{10}+2\theta^8+4\theta^6+8\theta^4+16\theta^2+32$, $q_2(\theta)=\theta^8+4\theta^6+12\theta^4+32\theta^2+80$, $q'_1(\theta)=\theta^8-2\theta^6+8\theta^2-16$ and $q'_2(\theta)=\theta^4-4\theta^2+4.$\\
 \indent Since the above $2$-integral basis is not in triangular form, we will find a new $2$-integral basis that is in triangular form. Clearly $\theta^{12}=m$ and $v_2(m)=6$ implies that $\frac{\theta^2}{2},$ $\frac{\theta^3}{2}$ and $\frac{\theta^4}{2^2}$ are algebraic integers. Since   $$\frac{q_1(\theta)}{2^6}-\frac{\theta^2q_1'(\theta)}{2^6}-\frac{q_2(\theta)}{2^4}=\frac{\theta^6+8}{2^4}-\frac{4\theta^6+80}{2^4}-\frac{12\theta^4+24\theta^2}{2^4},$$  we see that  $\frac{\theta^6+8}{2^4}$ is an algebraic integer. Hence in view of Proposition \ref{prop1}, the following set is a $2$-integral basis of $K.$ $$A_7=\left \{ 
 	1,\theta,\frac{\theta^2}{2},
 	\frac{\theta^3}{2}, \frac{\theta^4}{2^2}, \frac{\theta^5}{2^2}, \frac{\theta^6+8}{2^4}, \frac{\theta^7+8\theta}{2^4}, \frac{ q'_1(\theta)}{2^5},\frac{\theta q'_1(\theta)}{2^5},    
 	\frac{q_1(\theta)}{2^6}, \frac{\theta q_1(\theta)}{2^6} \right\}.$$\\
\textbf{Case A8:} $v_2(m)=6$, $m_2\equiv 3~(8)$. Arguing as in Case A7, we see that $f(x)$ is not $2$-regular. Here $\lambda$, $\psi_1(Y)$, $\psi_2(Y)$ and $N_1$ are same as in Case A7.\\
\indent Consider the first type $(x;-\frac{1}{2},Y+\bar{1}).$ Take $\Phi(x)=x^2+2x-2$ and $V_1$ be same as in Case A7, then $V_1(\Phi)=2$. The $V_1$-Newton polygon of 
\begin{equation}
\begin{split}	f(x)=\Phi^6-(12x-72)\Phi^5-(400x-1220)\Phi^4-(4512x-8832)\Phi^3-\\(23040x-31536)\Phi^2-(54976x-54656)\Phi-49920x+36544-m
\end{split}
\end{equation}
 has a single edge joining the points $(0,~17)$ and $(2,~12)$. The  residual polynomial associated to this edge is linear. Therefore $f(x)$ is $V_1$-regular and $N_2^1=2.$\\
\indent Consider the second type $(x;-\frac{1}{2},Y^2+Y+\bar{1})$. Take $\Phi(x)=x^4+2x^3-2x^2+4x+4$ and $V_2$ be same as in Case A7. For this type $V_2(\Phi)=4$. The $V_2$-Newton polygon of 
\begin{equation}
	\begin{split}
		f(x)=	\Phi^3(x)-(6x^3-18x^2+68x-224)\Phi^2(x)-(800x^3-1152x^2+2144x-3520)\Phi(x)\\-(13184x^3-12800x^2+11776x+17728+m)       
	\end{split}
\end{equation}
is the lower convex hull of the points $(0,~17)$, $(1,~15)$, $(2,~12)$ and $(3,~12).$ It has a single edge  joining the points $(0,~17)$ and $(2,~12)$. The  residual polynomial   associated to this edge is linear. Thus $f(x)$ is $V_2$-regular and   $N^2_2=2$.\\
\indent Hence by virtue of the Theorem \ref{ABC}, $v_2(\ind\theta)=N_1+\deg(\psi_1(Y))N_2^1+\deg(\psi_2(Y))N_2^2=30+2+4=36$ and the following set (not in triangular form) is a $2$-integral basis of $K.$
{\footnotesize $$
	\left \{ 
	\frac{q'_2(\theta)}{2^2},\frac{\theta q'_2(\theta)}{2^2},
	\frac{\theta^2 q'_2(\theta)}{2^3}, \frac{\theta^3 q'_2(\theta)}{2^3}, \frac{ q_2(\theta)}{2^4}, \frac{\theta q_2(\theta)}{2^4}, \frac{q'_1(\theta)}{2^5},\frac{\theta q'_1(\theta)}{2^5},    
	\frac{q_1(\theta)}{2^6}, \frac{\theta q_1(\theta)}{2^6}, \frac{\theta^2 q'_1(\theta)}{2^6}, \frac{\theta^3 q'_1(\theta)}{2^6} \right\},$$}
where 
$q_1(\theta)=\theta^{10}-2\theta^9+6\theta^8-16\theta^7+44\theta^6-120\theta^5+328\theta^4-896\theta^3+2448\theta^2-6688\theta+18272$, $q_2(\theta)=\theta^8-4\theta^7+16\theta^6-56\theta^5+188\theta^4-608\theta^3+1920\theta^2-5952\theta+18192$, $q'_1(\theta)=\theta^8-2\theta^7+6\theta^6-20\theta^5+56\theta^4-168\theta^3+504\theta^2-1488\theta+4432$ and $q'_2(\theta)=\theta^4-4\theta^3+16\theta^2-64\theta+228.$\\
 \indent In this case $$\frac{q_1(\theta)}{2^6}-\frac{\theta^2q_1'(\theta)}{2^6}=\frac{\theta^7-3\theta^6+12\theta^5-44\theta^4+148\theta^3-496\theta^2-1672\theta+4568}{2^4}$$ implies that $\frac{\theta^7+\theta^6-4\theta^5+4\theta^4+4\theta^3+8\theta-8}{2^4}$ is an algebraic integer.
   Take $r(\theta)=\theta^7+\theta^6-4\theta^5+4\theta^4+4\theta^3+8\theta-8.$ Keeping in mind that $1,\theta,\theta^2,\frac{\theta^3}{2},\frac{\theta^4}{2^2}$ are algebraic integers and $\frac{\theta r(\theta)-q_1'(\theta))-3r(\theta)}{2^4}=\frac{-12\theta^6+32\theta^5-32\theta^4+160\theta^3-496\theta^2+1456\theta-4416}{2^4}-\frac{\theta^6-4\theta^5+8\theta^4+4\theta^3+8}{2^4},$ we see that $\frac{\theta^6-4\theta^5+8\theta^4+4\theta^3-8}{2^4}=$ $\frac{s(\theta)}{2^4}$ (say) is an algebraic integer. Thus Proposition \ref{prop1} provide that the set  $$A_8=\left \{ 
1,\theta,\frac{\theta^2}{2},
\frac{\theta^3}{2}, \frac{\theta^4}{2^2}, \frac{\theta^5}{2^2}, \frac{s(\theta)}{2^4}, \frac{r(\theta)}{2^4}, \frac{ q'_1(\theta)}{2^5},\frac{\theta q'_1(\theta)}{2^5},    
\frac{q_1(\theta)}{2^6}, \frac{\theta q_1(\theta)}{2^6} \right\}$$ is a $2$-integral basis of $K.$\\
\textbf{Case A9:} $v_2(m)=6$, $m_2\equiv 7~(8)$. In this case $f(x)\equiv x^{12}~(2)$. Proceeding same as in Case A7, one can check that $f(x)$ is not $2$-regular and $N_1=30$. In this case $\lambda$, $\psi_1(Y)$, $\psi_2(Y)$, $V_1$ and $V_2$ will be same as in Case A7.\\
\indent For the first type  $(x;\lambda,\psi_1(Y))$, let $\Phi(x)=x^2+2x+2$.  The $V_1$-Newton polygon of
\begin{equation}
	\begin{split}	f(x)=\Phi^6(x)-(12x-48)\Phi^5(x)-(160x-20)\Phi^4(x)-(32x+448)\Phi^3(x)\\+(384x+432)\Phi^2(x)-(192x)\Phi(x)-64-m
	\end{split}
\end{equation}
 is the lower convex hull of the points $(0,~2v_2(64+m))$, $(1,~15)$ and $(2,~12)$. If $m_2\equiv7~(16)$, then it has   a single edge, say $S'$, joining the points $(0,~18)$ and $(2,~12)$ with a lattice point $(1,~15)$ lying on it. The  residual polynomial associated to  $S'$ is $Y^2+Y+\bar{1}\in \F_2[Y]$. If $m_2\equiv15~(16)$, then the $V_1$-Newton polygon has two edges of negative slope. The first edge is the line segment joining the point $(0,~2v_2(m+64))$ with $(1,~15)$ and the second edge is the one joining the point $(1,~15)$ to $(2,~12)$. The  residual polynomial attached to each edge is linear. Thus $f(x)$ is $V_1$-regular and  $N^1_2=3$.\\
\indent For the second type  $(x;\lambda,\psi_2(Y))$, take $\Phi(x)=x^4+2x^3+2x^2+4x+4$. The $V_2$-Newton polygon of
 \begin{equation}
 f(x)=\Phi^3(x)-(6x^3-6x^2+20x-32)\Phi^2(x)-(64x^3+32x^2+96x+128)\Phi(x)-64-m
 	\end{equation}
is the lower convex hull of the points $(0,~2v_2(64+m))$, $(1,~15)$, $(2,~12)$ and $(3,12)$.
  If $m_2\equiv7~(16)$, then it  has only one edge of negative slope. The edge is the line segment  joining the points $(0,~18)$ and $(2,~12)$ with a lattice point $(1,15)$ lying on it. The residual polynomial  associated to this edge  has no repeated roots.
  If $m_2\equiv15~(16)$, then the  $V_2$-Newton polygon of $f(x)$ has two edges of negative slope. The first edge joins $(0,~2v_2(64+m))$ with $(1,~15)$ and the second edge is from $(1,~15)$ to $(2,~12)$. The  residual polynomial  associated with each edge  is linear. Thus $f(x)$ is $V_2$-regular and $N_2^2=3$.\\
 \indent Hence in view of Theorem \ref{ABC}, we have $v_2(\ind\theta)=N_1+\deg\psi_1(Y) N_2^1+\deg\psi_2(Y)N_2^2=30+3+6=39$. Therefore the set  {\footnotesize $$
  	\left \{ 
  	\frac{q'_2(\theta)}{2^2},\frac{\theta q'_2(\theta)}{2^2},
  	\frac{\theta^2 q'_2(\theta)}{2^3}, \frac{\theta^3 q'_2(\theta)}{2^3}, \frac{ q_2(\theta)}{2^4}, \frac{\theta q_2(\theta)}{2^4}, \frac{q'_1(\theta)}{2^5},\frac{\theta q'_1(\theta)}{2^6},    
  	\frac{q_1(\theta)}{2^6}, \frac{\theta^2 q'_1(\theta)}{2^6}, \frac{\theta q_1(\theta)}{2^7}, \frac{\theta^3 q'_1(\theta)}{2^7} \right\},$$} where 
  $q_1(\theta)=\theta^{10}-2\theta^9+2\theta^8-4\theta^6+8\theta^5-8\theta^4+16\theta^2-32\theta+32$, $q_2(\theta)=\theta^8-4\theta^7+8\theta^6-8\theta^5-4\theta^4+32\theta^3-64\theta^2+64\theta+16$, $q'_1(\theta)=\theta^8-2\theta^7+2\theta^6-4\theta^5+8\theta^4-8\theta^3+8\theta^2-16\theta+16$ and $q'_2(\theta)=\theta^4-4\theta^3+8\theta^2-16\theta+36,$ is a $2$-integral basis of $K.$\\ 
  \indent One can easily check that $$\frac{q'_1(\theta)}{2^5}+\frac{\theta^3 q_1'(\theta)}{2^7}-\frac{\theta q_1(\theta)}{2^7}=\frac{\theta^7-2\theta^6+4\theta^4-8\theta^3-24\theta+16}{2^5}+\frac{\theta^2}{2}.$$ Let $q_3(\theta)=\theta^7-2\theta^6+4\theta^4-8\theta^3-24\theta+16,$ then we get $$\frac{q_1(\theta)}{2^6}-\frac{\theta^2q'_1(\theta)}{2^6}-\frac{q_3(\theta)}{2^4}=-\frac{\theta^6-4\theta^5+8\theta^4-12\theta^3-16\theta+8}{2^4}.$$ Set  $q_4(\theta)=\theta^6-4\theta^5+8\theta^4-12\theta^3-16\theta+8.$ As $\frac{\theta^2}{2}$ is an algebraic integer in $K,$ therefore  $\frac{q_3(\theta)}{2^5}$ and $\frac{q_4(\theta)}{2^4}$ are algebraic integers. Hence using Proposition \ref{prop1}, we see that $$A_9=\left \{1,\theta,\frac{\theta^2}{2}, \frac{\theta^3}{2}, \frac{\theta^4}{2^2}, \frac{\theta^5}{2^2}, \frac{q_4(\theta)}{2^4}, \frac{q_3(\theta)}{2^5}, \frac{q_1'(\theta)}{2^5}, \frac{\theta q_1'(\theta)}{2^6},\frac{q_1(\theta)}{2^6}, \frac{\theta q_1(\theta)}{2^7}\right\},$$ is a $2$-integral basis of $K.$\\
  
 \textbf{Case A10:} $v_2(m)=8$, $m_2\equiv1~(8)$. In this case $f(x)\equiv x^{12}~(2)$. The $2$-Newton polygon of $f(x)$  has  a single edge joining the points $(0,~8)$ and $(12,~0)$ with slope $\lambda=-\frac{2}{3}$. The residual polynomial attached to this edge  is $(Y+\bar{1})^4\in \F_2[Y]$. By using Lemma \ref{lemma1}, we have  $N_1=40$. Set $\psi(Y)=Y+\bar{1}$ and  $\Phi(x)=x^3+4$. In view of Equation \eqref{eq3}, we  define second order valuation $V$ attached to the data $(x;\lambda,\psi(Y))$ such that  $V(x)=2$,  $V(2)=3$ and  $V(\Phi)=6$. The $V$-Newton polygon of
 \begin{equation}\label{eq3.11}
 	f(x)=\Phi^4(x)-16\Phi^3(x)+96\Phi^2(x)-256\Phi(x)+256-m
  \end{equation} 
 is the lower convex hull of the points $(0,~3v_2(256-m))$, $(1,~30)$, $(2,~27)$, $(3,~30)$, $(4,~24)$. If $m_2\equiv 9~(16)$, then it has  two  edges of negative slope and % The first edge, say $S_{11}$, is the line segment joining the points $(0,~33)$ and $(2,~27)$ with a lattice point $(1,~30)$ lying on it. The second  edge, say $S_{12}$, is the one joining the points $(2,~27)$ with $(4,~24)$. 
 the  residual polynomial  corresponding to each edge is separable. % Now we have  $y_u=-\frac{2}{3}u+8$, where $10\le u\le 12$, $e=3, a_1=0, b_1=2, a_2=2, b_2=4, e_1=1, e_2=2, f_1=2, f_2=1, Y_1=30, Y_2=27, Y_3=\frac{51}{2}$ and  $Y_4=24$.\\
 If $m_2\equiv1~(16)$, then the $V$-Newton polygon of $f(x)$ has  three  edges of negative slope and %. The first edge, say $S_{21}$, is the line segment joining the points $(0,~3v_2(256-m))$ and $(1,~30)$, the second edge, say $S_{22}$, is the one joining the points $(1,~30)$ with $(2,~27)$ and the third edge, say $S_{23}$, joins $(2,~27)$ and $(4,~24)$.
  the  residual polynomial associated to each edge is linear. Thus $f(x)$ is $V$-regular and  $N_2=10$. Hence by virtue of the Theorem Theorem \ref{ABC}, $v_2(\ind\theta)=N_1+\deg\psi(Y) N_2=40+10=50$ %Now we have  $y_u=-\frac{2}{3}u+8$, where $10\le u\le 12$, $e=3, a_1=0, b_1=1, a_2=1, b_2=2, a_3=2, b_3=4, e_1=e_2=1, e_3=2, f_i=1$, for $i=1,2,3$  $Y_1=30, Y_2=27, Y_3=\frac{51}{2}$ and  $Y_4=24$. Therefore 
  and a $2$-integral basis of $K$ is $$
A_{10}=\left \{1, \theta,  \frac{\theta^2}{2},
\frac{q_3(\theta)}{2^2},\frac{\theta q_3(\theta)}{2^3},
\frac{\theta^2 q_3(\theta)}{2^3},  \frac{q_2(\theta)}{2^5},\frac{\theta q_2(\theta)}{2^5},
\frac{\theta^2 q_2(\theta)}{2^6}, \frac{q_1(\theta)}{2^8}, \frac{\theta q_1(\theta)}{2^8}, \frac{\theta^2q_1(\theta)}{2^9}\right\},$$
where $q_1(\theta)=\theta^9-4\theta^6+16\theta^3-64$, $q_2(\theta)=\theta^6-8\theta^3+48$ and $q_3(\theta)=\theta^3-12.$\\\\
\hspace{-0.2468cm}\textbf{Case A11:} $v_2(m)=8$, $m_2\equiv5~(8)$. Arguing  as in Case A10, one can check that $f(x)$ is not $2$-regular. Here $\lambda$, $\psi(Y)$, $\Phi(x)$ and $V$  are same as in Case A10. Using the $\Phi$-expansion of $f(x)$ given in \eqref{eq3.11}, one can easily check that the $V$-Newton polygon of $f(x)$ has a single edge of negative slope and the residual polynomial attached to this edge is $Y^2+Y+\bar{1}\in F_2[Y]$. So $f(x)$ is $V$-regular and $N_2=8$. Hence using  Lemma \ref{lemma1} and Theorem \ref{ABC}, we observe that $v_2(\ind\theta)=N_1+\deg\psi(Y)N_2=40+8=48$ and  a $2$-integral basis of $K$ is $$
A_{11}=\left \{1, \theta, \frac{\theta^2}{2},  
\frac{q_3(\theta)}{2^2},\frac{\theta q_3(\theta)}{2^3},
\frac{\theta^2 q_3(\theta)}{2^3},  \frac{q_2(\theta)}{2^5},\frac{\theta q_2(\theta)}{2^5},
\frac{\theta^2 q_2(\theta)}{2^6}, \frac{q_1(\theta)}{2^7}, \frac{\theta q_1(\theta)}{2^8}, \frac{\theta^2q_1(\theta)}{2^8}\right\},$$
where $q_1(\theta)=\theta^9-4\theta^6+16\theta^3-64$, $q_2(\theta)=\theta^6-8\theta^3+48$ and $q_3(\theta)=\theta^3-12$.\\\\
 \textbf{Case A12:} $v_2(m)=8$, $m_2\equiv 3~(8)$. Here $f(x)=x^{12}~(2)$. Proceeding same  as in Case A10, we observe  that $f(x)$ is not $2$-regular. In this case $\lambda$, $\psi(Y)$, $\Phi(x)$ and  $V$ will be same as in Case A10. Using \eqref{eq3.11}, we see that the $V$-Newton polygon of $f(x)$  has a single edge joining the points $(0,~27)$ and $(4,~24)$. The  residual polynomial  corresponding to this edge is linear. So $f(x)$ is $V$-regular and $N_2=3$. Thus by Lemma \ref{lemma1} and Theorem \ref{ABC},  $v_2(\ind\theta)=N_1+\mu N_2=40+3=43$ and a $2$-integral basis of $K$ is $$
 A_{12}=\left \{1, \theta, \frac{\theta^2}{2}, 
 \frac{q_3(\theta)}{2^2},\frac{\theta q_3(\theta)}{2^2},
 \frac{\theta^2 q_2(\theta)}{2^3}, \frac{q_2(\theta)}{2^4},\frac{\theta q_2(\theta)}{2^5},
 \frac{\theta^2 q_2(\theta)}{2^5},\frac{q_1(\theta)}{2^6}, \frac{\theta q_1(\theta)}{2^7}, \frac{\theta^2q_1(\theta)}{2^8} \right\},$$
 where $q_1(\theta)=\theta^9-4\theta^6+16\theta^3-64$, $q_2(\theta)=\theta^6-8\theta^3+48$ and $q_3(\theta)=\theta^3-12$.\\\\
 \noindent\textbf{Case A13:} $v_2(m)=10$, $m_2\equiv1~(4)$. Clearly $f(x)\equiv x^{12}~(2)$. The $2$-Newton polygon of $f(x)$  has a single edge  joining the points $(0,~10)$ and $(12,~0)$ with slope $\lambda=-\frac{5}{6}$. The residual polynomial  of $f(x)$ with respect to $2$ is $(Y+\bar{1})^2\in \F_2[Y]$. So $f(x)$ is not $2$-regular.
Let $\psi(Y)=Y+\bar{1}$. Take $\Phi(x)=x^6+32$, then $\Phi(x)$ works as the key polynomial attach to $\lambda$. The second order valuation $V$ attached to the data $(x;\lambda,\psi(Y))$ on $\Q_2(x)$ is defined as $V\big(\displaystyle\sum_{i\geq 0}^{}a_ix^i\big)=6\min_{i\ge0}\{v_2(a_i)+i\frac{5}{6}\}$. So $V(x)=5$, $V(2)=6$ and $V(\Phi)=30$.
Here  $f(x)=x^{12}-m=(\Phi(x)-32)^2-m=\Phi^2(x)-64\Phi(x)+1024-m.$
If $m_2\equiv1~(8)$, then the $V$-Newton polygon of $f(x)$  has two edges of negative slope and the  residual polynomial  attached to each edge is linear.
If $m_2\equiv5~(8)$, then the $V$-Newton polygon of $f(x)$  has a single edge of negative slope. The residual polynomial  corresponding to this edge is $Y^2+Y+\bar{1}\in \F_2[Y]$. Thus $f(x)$ is $V$-regular and $N_2=6$.
%Here  $y_u=-\frac{5}{6}u+10$, where $7\le u\le 12$, $e=6$, $a_1=0$, $b_1=2$, $e_1=1$, $f_1=2$,  $Y_1=66$,  $Y_2=60$ and $q_1(\theta)=\theta^6-32$.
 Hence in view of Lemma \ref{lemma1} and Theorem \ref{ABC} $v_2(\ind\theta)=N_1+\deg\psi(Y)N_2=50+6=56$ and $$
A_{13}=\left \{1, \theta, \frac{\theta^2}{2}, \frac{\theta^3}{2^2}, \frac{\theta^4}{2^3}, \frac{\theta^5}{2^4}, \frac{q_1(\theta)}{2^6}, \frac{\theta q_1(\theta)}{2^6}, \frac{\theta^2q_1(\theta)}{2^7}, \frac{\theta^3q_1(\theta)}{2^8}, \frac{\theta^4q_1(\theta)}{2^9}, \frac{\theta^5q_1(\theta)}{2^{10}}\right\},$$ where $q_1(\theta)=\theta^6-32,$ is a $2$-integral basis of $K.$\\\\
\textbf{Case A14:}  $v_2(m)=10$, $m_2\equiv 3~(8)$. Here $f(x)\equiv x^{12}~(2)$. Arguing as in Case A13, we see that $f(x)$ is not $2$-regular. In this case $\lambda$, $V$, $N_1$ and $\psi(Y)$  will be same as in Case A13. Take $\Phi(x)=x^6+8x^3-32$. The $V$-Newton polygon of  $f(x)=\Phi^2(x)+(128-16x^3)\Phi(x)-1024x^3+3072-m,$
  has a single edge joining the points $(0,~75)$ and $(2,~60)$. The  residual polynomial  associated to this edge is linear. Therefore $f(x)$ is $V$-regular and $N_2=7$. Hence by Theorem \ref{ABC}, $v_2(\ind\theta)=N_1+\mu N_2=50+7=57$ and a $2$-integral basis of $K$ is 
  $$A_{14}=\left \{1, \theta, \frac{\theta^2}{2}, \frac{\theta^3}{2^2}, \frac{\theta^4}{2^3}, \frac{\theta^5}{2^4}, \frac{q_1(\theta)}{2^6}, \frac{\theta q_1(\theta)}{2^{7}}, \frac{\theta^2q_1(\theta)}{2^{7}}, \frac{\theta^3q_1(\theta)}{2^{8}}, \frac{\theta^4q_1(\theta)}{2^{9}}, \frac{\theta^5q_1(\theta)}{2^{{10}}}\right\},$$
where $q_1(\theta)=\theta^6-8\theta^3+96$.\\\\
\noindent\textbf{Case A15:} $v_2(m)=10$, $m_2\equiv 7~(8)$. Here $f(x)\equiv x^{12}~(2)$. Proceeding same as in Case A13, we see that $f(x)$ is not $2$-regular. In this case $\lambda$, $\psi(Y)$ and $V$ will be same as in Case A13.  Take $\Phi(x)=x^6+8x^3+32$, then $V(\Phi)=30$. The $\Phi$-expansion of $f(x)$ is $f(x)=\Phi^2(x)-16x^3\Phi(x)-1024-m.$
If $m_2\equiv7~(16)$, then the $V$-Newton polygon  has a single edge joining the points $(0,~78)$ and $(2,~60)$ with a lattice point $(1,~69)$ lying on it. The residual polynomial  attached to this edge is $Y^2+Y+\bar{1}\in \F_2[Y]$.
If $m_2\equiv1~(16)$,  then the $V$-Newton polygon has two edges of negative slope and the  residual polynomial attached with each  edge is linear. Hence $f(x)$ is $V$-regular and $N_2=9$. Thus using Lemma \ref{lemma1} and Theorem \ref{ABC}, we see that $v_2(\ind\theta)=N_1+\mu N_2=50+9=59$ and  $$
A_{15}=\left \{1, \theta, \frac{\theta^2}{2}, \frac{\theta^3}{2^2}, \frac{\theta^4}{2^3}, \frac{\theta^5}{2^4}, \frac{q_1(\theta)}{2^6}, \frac{\theta q_1(\theta)}{2^{7}}, \frac{\theta^2q_1(\theta)}{2^{8}}, \frac{\theta^3q_1(\theta)}{2^{9}}, \frac{\theta^4q_1(\theta)}{2^{9}}, \frac{\theta^5q_1(\theta)}{2^{{10}}}\right\},$$
where $q_1(\theta)=\theta^6-8\theta^3+32,$ is a $2$-integral basis of $K.$\\\\
\noindent\textbf{Case B1:} $v_3(m)=3$,  $m_3\equiv \delta ~(9)$, where $\delta\in\{1,-1\}$. Here $f(x)\equiv x^{12}~(3)$. The $3$-Newton polygon of $f(x)$  has only one edge joining the points $(0,~3)$ and $(12,~0)$ having slope $-\frac{1}{4}$. The residual polynomial  attached with this edge is $(Y-\bar{\delta})^3\in \F_3[Y]$, where $\delta$=1 or -1 according as $m_3\equiv1~(3)$ or $m_3\equiv-1~(3)$. Let $\psi(Y)=Y-\bar{\delta}$. Take $\Phi(x)=x^4-3\delta,$ then it works as a key polynomial attached to $\lambda.$ In view of Equation \eqref{eq3}, we have a second order valuation $V$ on $\Q_3(x)$  defined as  $V\big(\displaystyle\sum_{i\geq 0}^{}a_ix^i\big)=4 \min_{i\ge0}\{v_3(a_i)+i\frac{1}{4}\}$ such that   $V(x)=1$, $V(3)=4$ and $V(\Phi)=4$. The $V$-Newton polygon of
\begin{equation}\label{eq3.15}
f(x)=\Phi^3(x)+9\delta\Phi^2(x)+27\Phi(x)+27\delta-m
\end{equation}
\begin{figure}[H]
	\centering
		\begin{tikzpicture}[scale=0.7]
			\draw[thick, ->] (-0.5, 0) -- (9, 0);
			\draw[thick, ->] (0, -0.5) -- (0, 8);
			\draw[thick] (0, 7) .. controls (2,4) and  (2,4) ..(2,4);
			\draw[thick] (2, 4) .. controls (7,2) and  (7,2) ..(7,2);
			
			\draw[dashed] (0,2) .. controls (7,2) and  (7,2) ..(7,2);

			\draw (0, 7) node[right]{$(0,4v_3(27\delta-m))$};
			\draw (-.35, 6.9) node[right]{$\bullet$};
			\draw (3.8,4.6) node[left]{$(1,16)$};
			\draw (7, 2) node[right]{$(3,12)$};
			
			\draw (0, 2) node[left]{$(0,12)$};
			\draw (7.3, 2) node[left]{$\bullet$};
			\draw (1.6, 2.5) node[right]{$\star$};
			\draw (1.6, 3) node[right]{$\star$};
			\draw (1.6, 3.5) node[right]{$\star$};
			\draw (1.6, 4.01) node[right]{$\star$};
			\draw (3.8, 3) node[right]{$\star$};
			\draw (3.8, 2.5) node[right]{$\star$};
			
		\end{tikzpicture}
		\caption{ $V$-Newton polygon of $f(x).$}
		\label{fig1.2}
\end{figure}
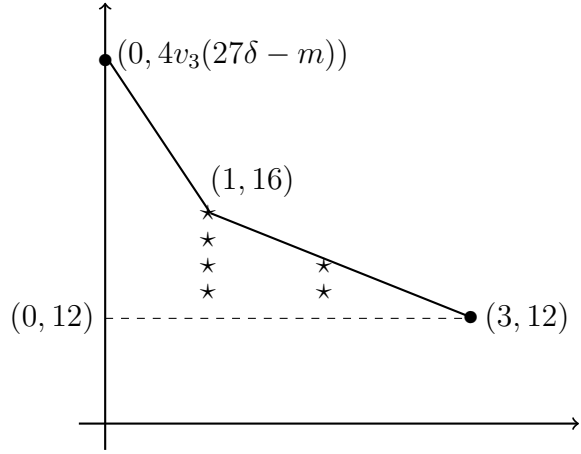
 \noindent has two edges of negative slope (See Figure \ref{fig1.2}). The first edge is the line segment joining the point $(0,~4v_3(27\delta-m))$ with $(1,~16)$ and the second edge is the one joining the point $(1,~16)$ to $(3,~12)$. The residual polynomial   associated to each edge has no repeated roots. Thus $f(x)$ is $V$-regular and $N_2=6$. According to the notation of Theorem \ref{ABC}, we have  $y_u=-\frac{1}{4}u+3$, where $9\le u\le 12$, $[a_1,b_1]=[0,1]$, $e_1=1$, $f_1=1$, $[a_2,b_2]=[1,3]$, $e_2=1$, $f_2=2$, $Y_j=-2j+18$, where $1\le j \le 3$. Hence using Lemma \ref{lemma1} and Theorem \ref{ABC}, we see that $v_3(\ind\theta)=N_1+\deg\psi(Y) N_2=18=12+6=18$ and a $3$-integral basis of $K$ is $$
B_1=\left \{1, \theta, \theta^2, \theta^3, \frac{q_2(\theta)}{3}, \frac{\theta q_2(\theta)}{3}, \frac{\theta^2q_2(\theta)}{3^2}, \frac{\theta^3q_2(\theta)}{3^2}, \frac{q_1(\theta)}{3^3}, \frac{\theta q_1(\theta)}{3^{3}}, \frac{\theta^2q_1(\theta)}{3^3}, \frac{\theta^3q_1(\theta)}{3^{3}}\right\},$$
where $q_1(\theta)=\theta^8+3\delta\theta^4+9$ and $q_2(\theta)=\theta^4+6\delta$.\\\\
\noindent\textbf{Case B2:} $v_3(m)=3$,  $m_3\in\{4\delta,7\delta\}~(9),~ \delta\in\{1,-1\}$. In this case $f(x)\equiv x^{12}~(3)$. Arguing as in Case B1, we see that $f(x)$ is not $3$-regular. In this case $\lambda$, $\psi(Y)$, $\phi(x)$ and  $V$ will be same as in Case B1. Keeping in mind \eqref{eq3.15}, we observe that the $V$-Newton polygon of $f(x)$   has a single edge joining the points  $(0,~16)$ and $(3,~12)$. The residual polynomial associated with this edge is linear. So $f(x)$ is $V$-regular and $N_2=3$. Therefore Lemma \ref{lemma1} and Theorem \ref{ABC} implies that $v_3(\ind\theta)=N_1+\deg\psi(Y) N_2=12+3=15$ and the set $$
B_2=\left \{1, \theta, \theta^2, \theta^3, \frac{q_2(\theta)}{3}, \frac{\theta q_2(\theta)}{3}, \frac{\theta^2q_2(\theta)}{3}, \frac{\theta^3q_2(\theta)}{3^2}, \frac{q_1(\theta)}{3^2}, \frac{\theta q_1(\theta)}{3^{2}}, \frac{\theta^2q_1(\theta)}{3^3}, \frac{\theta^3q_1(\theta)}{3^{3}}\right\},$$
where $q_1(\theta)=\theta^8+3\delta\theta^4+9$ and $q_2(\theta)=\theta^4+6\delta$, is a $3$-integral basis of $K.$\\\\
\noindent\textbf{Case B3:} $v_3(m)=6$, $m_3\equiv 1~(9)$. Here $f(x)\equiv x^{12}~(3)$. The $3$-Newton polygon of $f(x)$ has a single edge joining the points $(0,~6)$ and (12,~0) with slope $-\frac{1}{2}$. The residual polynomial  associated to this edge is $Y^6-\bar{1}=(Y-\bar{1})^3(Y+\bar{1})^3\in \F_3[Y]$, which has a repeated root. Therefore $f(x)$ is not $3$-regular. At this step, we have two types of data $(x;-\frac{1}{2},Y-\bar{1})$ and $(x;-\frac{1}{2}, Y+\bar{1}).$ Set $\psi(Y)=Y+\bar{\delta}$, where $\delta\in \{1,-1\}$.   Take $\Phi(x)=x^2+3\delta$. In view of Equation \eqref{eq3}, we can define the second order valuation $V$ on $\Q_3(x)$  as $V\big(\displaystyle\sum_{i\geq 0}^{}a_ix^i\big)=2 \min_{i\ge0}\{v_3(a_i)+i\frac{1}{2}\}$ so that $V(\Phi)=2$, $V(x)=1$ and $V(3)=2$. The $\Phi$-expansion of $f(x)$ is given by 
\begin{equation}\label{EQ}
f(x)=\Phi^6(x)-18\delta\Phi^5(x)+135\Phi^4(x)-540\delta\Phi^3(x)+1215\Phi^2(x)-1458\delta\Phi(x)+729-m.
\end{equation}
The $V$-Newton polygon of $f(x)$  being the lower convex hull of the points $(0,~2v_3(729-m)),$ $(1,~14),$ $(2,~14),$ $(3,~12),$ $(4,~14),$ $(5,~14),$ $(6,~12)$ has two edges of negative slope. The first edge is the line segment joining the point $(0,~2v_3(729-m))$ with $(1,~14)$ and the second edge is the one joining the point $(1,~14)$ to $(3,~12)$. The residual polynomial   associated with each edge is separable. Therefore $f(x)$ is $V$-regular and $N_2=3$.
Hence using Lemma \ref{lemma1} and Theorem \ref{ABC}, we see that $v_3(\ind\theta)=N_1+\deg(Y-1)N_2+\deg(Y+1)N_2=30+3+3=36$ and a $3$-integral basis of $K$ is 
{\small $\left \{\frac{q_3(\theta)}{3^3}, \frac{\theta q_3(\theta)}{3^3}, \frac{q'_3(\theta)}{3^3}, \frac{\theta q'_3(\theta)}{3^3}, \frac{q_2(\theta)}{3^4}, \frac{ q'_2(\theta)}{3^4},  \frac{\theta q_2(\theta)}{3^5}, \frac{\theta q'_2(\theta)}{3^5}, \frac{q_1(\theta)}{3^6}, \frac{\theta q_1(\theta)}{3^6}, \frac{ q'_1(\theta)}{3^6}, \frac{\theta q'_1(\theta)}{3^6} \right\},$}
where $q_1(\theta)=\theta^{10}+3\theta^8+9\theta^6+27\theta^4+81\theta^2+243$, $q_2(\theta)=\theta^8+6\theta^6+27\theta^4+108\theta^2+405$ and  $q'_1(\theta)=\theta^{10}-3\theta^8+9\theta^6-27\theta^4+81\theta^2-243,$ $q'_2(\theta)=\theta^8-6\theta^6+27\theta^4-108\theta^2+405$, $q_3(\theta)=\theta^6+9\theta^4+54\theta^2+270$ and $q'_3(\theta)=\theta^6-9\theta^4+54\theta^2-270.$\\
\indent Since $\theta^{12}=m$ and $v_3(m)=6,$ we conclude that $1,\theta,\frac{\theta^2}{3},\frac{\theta^3}{3}, \frac{\theta^4}{3^2}$ are algebraic integers. As $\frac{q_3(\theta)}{3^3}, \frac{\theta q_3(\theta)}{3^3}$ and $\frac{q_2(\theta)}{3^4}$ are algebraic integers, so  $\frac{\theta^6}{3^3}, \frac{\theta^7}{3^3}$ and $\frac{\theta^8}{3^4}$ are algebraic integers. Thus  $\frac{\theta^8-18\theta^4-162}{3^5}=\frac{q_1'(\theta)}{3^6}-\frac{q_1(\theta)}{3^6}+\frac{\theta^8}{3^4}$ is an algebraic integer. Let $z_2(\theta)=\theta^8-18\theta^4-162,$ then $\frac{\theta z_2(\theta)-\theta q_2(\theta)}{3^5}=-\frac{3\theta^7+27\theta^3}{3^4}+\frac{\theta^7-15\theta^5-9\theta^3-27\theta}{3^4}-2\theta,$ shows that $\frac{\theta^7-15\theta^5-9\theta^3-27\theta}{3^4}=\frac{z_1(\theta)}{3^4}$ (say) is an algebraic integer. Also $$\frac{\theta q_2(\theta)-\theta q'_2(\theta)}{3^5}-\frac{z_1(\theta)}{3^4}=\frac{6\theta^5}{3^3}+\frac{3\theta^7+81\theta^3}{3^4}-\frac{\theta^5-9\theta}{3^3}$$ implies that  $\frac{\theta^5-9\theta}{3^3}$ is an algebraic integer. Hence by using Proposition \ref{prop1}, the set 
 $$B_3=\left\{1,\theta, \frac{\theta^2}{3}, \frac{\theta^3}{3},\frac{\theta^4}{3^2}, \frac{\theta^5-9\theta}{3^3}, \frac{ \theta^6-9\theta^2}{3^3}, \frac{ z_1(\theta)}{3^4},\frac{z_2(\theta)}{3^5}, \frac{\theta z_2(\theta)}{3^5}, \frac{q_1(\theta)}{3^6}, \frac{\theta q_1(\theta)}{3^6}\right\},$$ is a $3$-integral basis  of $K.$ \\
 
\noindent\textbf{Case B4:} $v_3(m)=6$, $m_3\in\{4,7\}~(9)$. Now $f(x)\equiv x^{12}~(3)$. Proceeding same as in Case B3, it is easy to check that $f(x)$ is not $3$-regular. Here $\lambda$, $\psi(Y)$, $\Phi(x)$ and $V$ will be same as in Case B3. Using the $\Phi$-expansion of $f(x)$ given in \eqref{EQ}, the $V$-Newton polygon of $f(x)$  being the lower convex hull of the points of the set $T=\{(0,~14),(1,~14),(2,~14),(3,~12),(4,~14),(5,~14),(6,~12)\}$ has a single edge of negative slope. The edge is the line segment joining the points $(0,~14)$ and $(3,~12).$ The  residual polynomial associated with this edge is linear. Therefore $f(x)$ is $V$-regular and  $N_2=1.$ Hence  $v_3(\ind\theta)=N_1+ 2N_2=30+2=32$  and a $3$-integral basis (not in triangular form) of $K$ is 
 {\small $B_4=\left \{\frac{q_3(\theta)}{3^3}, \frac{\theta q_3(\theta)}{3^3}, \frac{ q'_3(\theta)}{3^3}, \frac{\theta q'_3(\theta)}{3^3}, \frac{q_2(\theta)}{3^4}, \frac{ \theta q_2(\theta)}{3^4}, \frac{q'_2(\theta)}{3^4}, \frac{\theta q'_2(\theta)}{3^4}, \frac{q_1(\theta)}{3^5}, \frac{q'_1(\theta)}{3^5}, \frac{\theta q_1(\theta)}{3^6}, \frac{\theta q'_1(\theta)}{3^6} \right\},$}
where $q_1(\theta)=\theta^{10}+3\theta^8+9\theta^6+27\theta^4+81\theta^2+243$, $q_2(\theta)=\theta^8+6\theta^6+27\theta^4+108\theta^2+405$ and  $q'_1(\theta)=\theta^{10}-3\theta^8+9\theta^6-27\theta^4+81\theta^2-243$ and $q'_2(\theta)=\theta^8-6\theta^6+27\theta^4-108\theta^2+405$, $q_3(\theta)=\theta^6+9\theta^4+54\theta^2+270$ and $q'_3(\theta)=\theta^6-9\theta^4+54\theta^2-270.$\\
\indent In this case $\theta q_1(\theta)-\theta q_1'(\theta)=6\theta^9+54\theta^5+486\theta$ implies that $\frac{\theta^9-18\theta^5-162\theta}{3^5}$ is an algebraic integer. Hence a triangular $3$-integral basis of $K$ is $$ B_4=\left\{1,\theta, \frac{\theta^2}{3}, \frac{\theta^2}{3},\frac{\theta^4}{3^2}, \frac{\theta^5}{3^2}, \frac{\theta^6}{3^3}, \frac{\theta^7}{3^3},\frac{\theta^8}{3^4}, \frac{\theta^9-18\theta^5-162\theta}{3^5}, \frac{q_1(\theta)}{3^5}, \frac{\theta q_1(\theta)}{3^6}\right\}.$$\\\\
\textbf{Case B5:} $v_3(m)=6$, $m_3\equiv -1\mod9$. Note that $f(x)\equiv x^{12}\mod3$. Here the $3$-Newton polygon of $f(x)$  has a single edge joining the points $(0,~6)$ and (12,~0) with slope $-\frac{1}{2}$. The  residual polynomial of $f(x)$ associated to the edge is $Y^6+\bar{1}=(Y^2+\bar{1})^3$. Kepping in mind Lemma \ref{lemma1}, we have $N_1=30$. Set $\psi(Y)=Y^2+\bar{1}$. 
Take  $\Phi(x)=x^4+9$. In view of Equation \eqref{eq3}, we have $V\big(\displaystyle\sum_{i\geq 0}^{}a_ix^i\big)=2 \min_{i\ge0}\{v_3(a_i)+i\frac{1}{2}\}$ with
 $V(\Phi)=4$, $V(x)=1$ and $V(3)=2$. The $V$-Newton polygon of \begin{equation}\label{eq3.17}
 	f(x)=\Phi^3(x)-27\Phi^2(x)+243\Phi(x)-729-m
 \end{equation}
 has two edges of negative slope. The first edge is the line segment joining the point $(0,~2v_3(m+729))$ with $(1,~14)$ and the second edge is the one joining the point  $(1,~14)$ to $(3,~12)$. The residual polynomial  associated to each edge is separable.  Hence $v_3(\ind\theta)=N_1+\deg\psi(Y)N_2=36$ and a $3$-integral basis of $K$ is $$
B	_5=\left \{1, \theta, \frac{\theta^2}{3}, \frac{\theta^3}{3}, \frac{q_2(\theta)}{3^2}, \frac{\theta q_2(\theta)}{3^3}, \frac{\theta^2 q_2(\theta)}{3^3}, \frac{\theta^3 q_2(\theta)}{3^{4}}, \frac{q_1(\theta)}{3^{5}}, \frac{\theta q_1(\theta)}{3^{5}}, \frac{\theta^2 q_1(\theta)}{3^{6}}, \frac{\theta^3q_1(\theta)}{3^{{6}}}\right\},$$
where $q_1(\theta)=\theta^8-9\theta^4+81$ and $q_2=\theta^4-18$.\\\\
\textbf{Case B6:} $v_2(m)=6$, $m_3\in\{-4,-7\}~(9)$. Proceeding same as in Case B5, one can easily check that $f(x)$ is not $2$-regular. Let $\lambda$, $\psi(Y)$, $\Phi(x)$ and $V$ be same as in Case B5. By Equation \eqref{eq3.17}, we observe that the $V$-Newton polygon of $f(x)$  has a single edge of negative slope. The edge is the line segment joining the points $(0,~14)$ and $(3,~12)$. The residual polynomial  associated with the edge is linear. Therefore $v_3(\ind\theta)=N_1+\deg\psi(Y) N_2=32$ and a $3$-integral basis of $K$ is $$
B_6=\left \{1, \theta, \frac{\theta^2}{3}, \frac{\theta^3}{3}, \frac{q_2(\theta)}{3^2}, \frac{\theta q_2(\theta)}{3^2}, \frac{\theta^2 q_2(\theta)}{3^3}, \frac{\theta^3 q_2(\theta)}{3^3}, \frac{q_1(\theta)}{3^4}, \frac{\theta q_1(\theta)}{3^5}, \frac{\theta^2 q_1(\theta)}{3^5}, \frac{\theta^3q_1(\theta)}{3^6}\right\},$$
where $q_1(\theta)=\theta^8-9\theta^4+81$ and $q_2=\theta^4-18.$\\\\
\textbf{Case B7:} $v_3(m)=9$, $m_3\equiv\delta~(9)$, where $\delta\in\{1,-1\}$. The $3$-Newton polygon of $f(x)$ has a single edge of negative slope which is a line segment joining the points $(0,~9)$ and $(12,~0)$ having slope $-\frac{3}{4}$. The residual polynomial   is $(Y-\bar{\delta})^3$, where $\delta$ is $1$ or $-1$ according as $m_3\equiv1\mod3$ or $m_3\equiv -1\mod3$. Take $\psi(Y)=Y-\bar{\delta}$. Take  $\Phi(x)=x^4-27\delta$. The second order valuation on $\Q_3(x)$ is given by $V\big(\displaystyle\sum_{i\geq 0}^{}a_ix^i\big)=4 \min_{i\ge0}\{v_3(a_i)+i\frac{3}{4}\}$. Then $V(\Phi)=12$, $V(x)=3$ and $V(3)=4$. The $V$-Newton polygon of
\begin{equation}\label{eq3.18}
f(x)=\Phi^3(x)+81\delta\Phi^2(x)+2187\Phi(x)+19683\delta-m
\end{equation} 
 has   two edges of negative slope. The first edge is the line segment joining the point $(0,~4v_3(-m+19683\delta))$ with $(1,~40)$ and the second edge is the one joining the point $(1,~40)$ to $(3,~36)$. The residual polynomial attached to each edge is separable. In this case  $N_1=45$. Hence $v_3(\ind\theta)=N_1+\deg\psi(Y) N_2=51$  and a $3$-integral basis of $K$ is $$
B_7=\left \{1, \theta, \frac{\theta^2}{3}, \frac{\theta^3}{3^2}, \frac{q_2(\theta)}{3^3}, \frac{\theta q_2(\theta)}{3^4}, \frac{\theta^2q_2(\theta)}{3^5}, \frac{\theta^3q_2(\theta)}{3^5}, \frac{q_1(\theta)}{3^7}, \frac{\theta q_1(\theta)}{3^{7}}, \frac{\theta^2q_1(\theta)}{3^8}, \frac{\theta^3q_1(\theta)}{3^{9}}\right\},$$
where $q_1(\theta)=\theta^8+27\delta\theta^4+729$ and $q_2(\theta)=\theta^4+54\delta.$ \\\\
\textbf{Case B8:} $v_3(m)=9$, $m_3\in\{\pm4,\pm7\}~(9)$. Clearly $f(x)\equiv x^{12}~(3)$. With the argument given in  Case B7, we conclude that $f(x)$ is not $3$-regular. In this case $\lambda$, $\psi(Y)$, $\Phi(x)$ and $V$ will be same as in Case B7. Keeping in mind \eqref{eq3.18}, the $V$-Newton polygon of $f(x)$  has a single edge  joining the points $(0,~40)$ and $(3,~36)$. The residual polynomial  associated with the edge is linear. Hence $v_3(\ind\theta)=N_1+\deg\psi(Y) N_2=48$ and a $3$-integral basis of $K$ is $$
B_8=\left \{1, \theta, \frac{\theta^2}{3}, \frac{\theta^3}{3^2}, \frac{q_2(\theta)}{3^3}, \frac{\theta q_2(\theta)}{3^{4}}, \frac{\theta^2q_2(\theta)}{3^{4}}, \frac{\theta^3q_2(\theta)}{3^{5}}, \frac{q_1(\theta)}{3^6}, \frac{\theta q_1(\theta)}{3^7}, \frac{\theta^2q_1(\theta)}{3^8}, \frac{\theta^3q_1(\theta)}{3^8} \right\},$$
where $q_1(\theta)=\theta^8+27\delta\theta^4+729$ and $q_2(\theta)=\theta^4+54\delta.$\\
This completes the proof of the theorem.
 \end{proof}

\section{Proof of Theorem \ref{Th1.3} }
\begin{proof} \textbf{(1)}  In this case, the $p$-Newton polygon of $f(x)$ has a single edge joining the points $(0,~v_p(m))$ and $(12,~0).$ As $p\nmid v_p(m)$, so  the residual polynomial associated to this edge is separable. In view of Theorem \ref{SK}, we have $v_p(\ind\theta)=\frac{1}{2}[11(v_p(m)-1)+(d-1)].$ Write $m=p^{v_p(m)}m_1,$ where $p\nmid m_1,$ then $\theta^{12}=p^{v_p(m)}m_1.$ This implies that  $\frac{\theta^{i}}{p^{\lfloor{\frac{iv_p(m)}{12}\rfloor}}}$ is a root of the polynomial $x^{12}-m^i_1p^{iv_p(m)-12\lfloor{\frac{iv_p(m)}{12}\rfloor}}$ having integer coefficient. Therefore $\frac{\theta^{i}}{p^{\lfloor{\frac{iv_p(m)}{12}\rfloor}}}$ is an algebraic integer in $K.$ Thus Proposition \ref{prop1}, provides that 
	$S=\left\{\frac{\theta^{i}}{p^{\lfloor{\frac{iv_p(m)}{12}\rfloor}}} |~ 1  \le i\le 12\right\}$ is a $p$-integral basis of $K$.\\
	\textbf{(2)}  $p=3$ and $p\nmid m.$  Then $f(x)\equiv(x-1)^3(x+1)^3(x^2+1)^3 (3)$ or  $f(x)\equiv (x^2-x-1)^3(x^2+x-1)^3 (3),$ according as $m\equiv -1~(3)$ or $m\equiv 1~(3).$ Take $\phi_1(x)=x^2+1$, $\phi_2(x)=x-\delta,$ and $\phi_3(x)=x^2-\delta x-1.$ where $\delta\in\{1,-1\}.$ For each $i=1,2,3,$ the $\phi_i$-expansion of $f(x)$ is given by 
		\begin{equation}
			f(x)=\phi_1^6(x)-6\phi_1^5(x)+15\phi_1^4(x)-17\phi_1^3(x)+12\phi_1^2(x)-6\phi_1(x)+1-m,
		\end{equation}
		\begin{equation}\label{eqqqq}
			\begin{split}
			f(x)=\phi_2^{12}(x)+12\delta\phi_2^{11}(x)+66\phi_2^{10}(x)+220\delta\phi_2^{9}(x)+495\phi_2^{8}(x)+ 792\delta\phi_2^{7}(x)+\hspace*{0.3in}\\924\phi_2^{6}(x)+ 792\delta\phi_2^{5}(x)+495\phi_2^{4}(x)+220\delta\phi_2^{3}(x)+66\phi_2^{2}(x)+12\delta\phi_2(x)+1-m, 
			\end{split}
              \end{equation}
		\begin{equation}
			\begin{split}
				f(x)=\phi_3^6(x)+(6\delta x+21)\phi_3^5(x)+(65\delta x+125)\phi_3^4(x)+(256\delta x+338)\phi_3^3(x)\\+(474\delta x+468)\phi_3^2(x)
	+(42\delta x+324)\phi_3(x)+89-m+144\delta x.
\end{split}
\end{equation}
 \indent Let $m^2 \equiv 1~(9),$ then for each $i=1,2,3,$  the $\phi_i$-Newton polygon of $f(x)$ has two  edges of negative slope. The first edge is the line segment joining the point $(0,~1)$ with $(1,~1)$ and the second  edge is the one joining the point  $(1,~1)$ to $(3,~0).$ The residual polynomial associated to each edge is linear. In view of Lemma \ref{lemma1} and Theorem \ref{ABC},  $v_3(\ind\theta)=4.$
Write $\theta^{12}-m\equiv (\theta^4-m)(\theta^8+m\theta^4+m^2)(3)$ and take $\xi=m^2\theta^8+m\theta^4+1,$ then $(m\theta^4-1)\xi=m^4-1.$ Therefore we see that $(m\theta^4)^3\xi^3=((m\theta^4-1+1)\xi)^3=(m^4-1+\xi)^3.$ By virtue of binomial theorem, we get $(m^4-1)\xi^3+3\xi^2(m^4-1)+3\xi(m^2-1)^2+(m^4-1)^3=0.$ As $m^4\ne 1,$ so $\xi$ satisfies the equation $$x^3+x^2+3(m^4-1)x+(m^4-1)^2 = 0$$ having integer coefficient. Since $m^2\equiv 1~(9),$   $\frac{\xi}{3}$ is an algebraic integer in $K.$ Therefore $\xi-\frac{(m^4-1)\theta^8}{3}=\frac{\theta^8+m\theta^4+1}{3}$ is an algebraic integer in $K.$ Let $h(\theta)=\theta^8+m\theta^4+1,$ then Proposition \ref{prop1} implies that the set $\{1,\theta,\theta^2,\cdots,\theta^7, \frac{h(\theta)}{3}, \frac{\theta h(\theta)}{3}, \frac{\theta^2 h(\theta)}{3},\frac{\theta^3h(\theta)}{3}\},$ is a $3$-integral basis of $K.$\\
\indent Let $m^2\not \equiv 1~(9),$ then for each $i=1,2,3,$ $\phi_i$-Newton polygon of $f(x)$ has a single edge  joining the points $(0,~1)$ and $(3,~0).$ The residual polynomial associated to this edge is linear. In view of Lemma \ref{lemma1} and Theorem \ref{ABC},  $v_3(\ind\theta)=0$ and therefore the set $\{1,\theta,\theta^2,\cdots,\theta^{11}\}$ is a $3$-integral basis of $K.$\\
\textbf{(3)} $p=2$ and $p\nmid m,$ then $f(x)\equiv(x+1)^4(x^2+x+1)^4(2).$ Let $\phi_1(x)=x^2+x+1$ and  $\phi_2(x)=x+1$.  The $\phi_1$-expansion of $f(x)$ is given by  \begin{equation}
	\begin{split}
	f(x)=\phi_1^6(x)+(9-6x)\phi_1^5(x)-(5x+25)\phi_1^4(x)+(24x+18)\phi_1^3(x)-18x\phi_1^2(x)\\+(4x-4)\phi_1(x)+1-m.
\end{split}
\end{equation}
 Keeping in mind the above expansion and expansion given in \ref{eqqqq}, one can verify that for each $i=1,2,$ the $\phi_i$- Newton polygon of $f(x)$ is the lower convex hull of the points $(0,~v_2(1-m))$, $(1,~2)$, $(2,~1)$, $(3,~i)$ and $(4,~0).$\\
\indent If $m\equiv 1~(8),$ then for each $i=1,2,$ the $\phi_i$-Newton polygon of $f(x)$ has either $2$ or $3$ edges of negative slope. The residual polynomial associated to each edge is separable. In virtue of the Theorem \ref{ABC}, we have $v_2(\ind\theta)=3\deg\phi_1(x)+3\deg\phi_2(x)=9.$ Clearly $\theta^{12}-m\equiv(\theta^3-1)(\theta^9+\theta^6+\theta^3+1)~(2).$ Set $g(\theta)=\theta^9+\theta^6+\theta^3+1.$ Using the same procedure as in the previous case, we find  that $\frac{g(\theta)}{4}$ is a root of the polynomial $$x^4-x^3-\frac{3(m-1)}{8}x^2+\frac{(m-1)^2}{16}x+\frac{(m-1)^3}{256}.$$ Hence in view of Proposition \ref{prop1}, \{$1$, $\theta$, $\theta^2$, $\theta^3$ ,$\theta^4$, $\theta^5$, $\frac{\theta^6-1}{2}$, $\frac{\theta(\theta^6-1)}{2}$, $\frac{\theta^2(\theta^6-1)}{2}$, $\frac{g(\theta)}{2^2}$, $\frac{\theta g(\theta)}{2^2}$, $\frac{\theta^2g(\theta)}{2^2}$\}
 is a $2$-integral basis of $K$\\
\indent If $m\equiv 5~(8),$ then for each $i=1,2,$ the $\phi_i$-Newton polygon of $f(x)$ has two edges of negative slope. The residual polynomial associated to each edge is linear. Using Theorem \ref{ABC}, we have $v_2(\ind\theta)=2\deg\phi_1(x)+2\deg\phi_2(x)=6.$  Write $\theta^{12}=m,$ then $(\theta^6-1+1)^2=m$. So $\frac{\theta^6-1}{2}$ is a root of the polynomial $x^2+x+\frac{1-m}{4}.$  Hence by virtue of the Proposition \ref{prop1}, \{$1$, $\theta$, $\theta^2$, $\theta^3$ ,$\theta^4$, $\theta^5$, $\frac{\theta^6-1}{2}$, $\frac{\theta(\theta^6-1)}{2}$, $\frac{\theta^2(\theta^6-1)}{2}$, $\frac{\theta^3(\theta^6-1)}{2}$, $\frac{\theta^4(\theta^6-1)}{2}$, $\frac{\theta^5(\theta^6-1)}{2}$\} is a $2$-integral basis of $K.$\\
\indent If $m\in\{3,7\}~(8),$ then for each $i=1,2,$ the $\phi_i$-Newton polygon of $f(x)$ has a single edge of negative slope. The residual polynomial associated to this edge is linear. Using Lemma \ref{lemma1} and Theorem \ref{ABC}, we have $v_2(\ind\theta)=0$. Therefore $\{1, \theta,\cdots, \theta^{11}\}$ is a $2$-integral basis of $K.$\\
\textbf{(4)} $p\nmid 12m.$ In this case $p\nmid D_f ,$ therefore $v_p(\ind\theta)=0.$ Hence $\{1,\theta,\cdots,\theta^{11}\}$ is a $p$-integral basis of $K.$
\end{proof}

\end{document}